# NUMERICAL APPROXIMATIONS FOR A THREE COMPONENTS CAHN-HILLIARD PHASE-FIELD MODEL BASED ON THE INVARIANT ENERGY QUADRATIZATION METHOD

XIAOFENG YANG*, JIA ZHAO†, QI WANG‡, AND JIE SHEN§

**Abstract.** How to develop efficient numerical schemes while preserving the energy stability at the discrete level is a challenging issue for the three component Cahn-Hilliard phase-field model. In this paper, we develop first and second order temporal approximation schemes based on the "Invariant Energy Quadratization" approach, where all nonlinear terms are treated semi-explicitly. Consequently, the resulting numerical schemes lead to a well-posed linear system with the symmetric positive definite operator to be solved at each time step. We rigorously prove that the proposed schemes are unconditionally energy stable. Various 2D and 3D numerical simulations are presented to demonstrate the stability and the accuracy of the schemes.

**1. Introduction.** The phase-field (i.e. the diffuse-interface) method is a robust modeling approach to simulate the free interface problem (cf. [1,9,11,19,22,25,26,26,27,31,32,35,56,60,61,65,67,68], and the references therein). Its essential idea is to use one (or more) continuous phase-field variable(s) to describe different phases in the multi-components system (usually immiscible), and represents the interfaces by thin, smooth transition layers. The constitutive equations for those phase-field variables can be derived from the energetic variational formalism, the governing system of equations is thereby well-posed and thermo-dynamically consistent. Consequently, one can carry out mathematical analyses to obtain the existence and uniqueness of solutions in appropriate functional spaces, and to develop convergent numerical schemes.

For an immiscible, two phase system, the commonly used free energy for the system consists of (i) a double-well bulk part which promotes either of the two phases in the bulk, yielding a hydrophobic contribution to the free energy; and (ii) a conformational capillary entropic term that promotes hydrophilic property in the multiphase material system. The competition between the hydrophilic and hydrophobic part in the free energy enforces the coexistence of two distinctive phases in the immiscible two-phase system. The corresponding binary system can be modeled either by the Allen–Cahn equation (second order) or the Cahn–Hilliard equation (fourth order). For both of these models, there have been many theoretical analysis, algorithm developments and numerical simulations available in literatures (cf. [1, 22, 25, 32]).

The generalization from the two-phase system to multi-phases have been studied by many works (cf. [2, 3, 5–7, 14, 15, 20, 21, 28–30, 38]). Specifically, for the system with three components, the general framework is to adopt three independent phase variables $(c_1, c_2, c_3)$ while imposing a hyperplane link condition among the three variables ($c_1 + c_2 + c_3 = 1$). The bulk part of the free energy is simply a summation of the original double-well energy for each phase variable [6, 7, 27]. Moreover, in order to ensure the boundedness (from below) of the free energy, especially for the *total spreading* case

*Corresponding author, Department of Mathematics, University of South Carolina, Columbia, SC 29208, USA. Email: xfyang@math.sc.edu. This author's research is partially supported by the U.S. National Science Foundation under grant numbers DMS-1200487 and DMS-1418898.

†Department of Mathematics, University of North Carolina at Chapel Hill, Chapel Hill, NC, 27599, USA. Email: zhaojia@email.unc.edu.

‡Department of Mathematics, University of South Carolina, Columbia, SC 29208; Beijing Computational Science Research Center, Beijing, China; School of Mathematics, Nankai University, Tianjin, China; Email: qwang@math.sc.edu. This author's research is partially supported by grants of the U.S. National Science Foundation under grant number DMS-1200487 and DMS-1517347, and a grant from the National Science Foundation of China under the grant number NSFC-11571032.

§Department of Mathematics, Purdue University, West Lafayette, IN 47906; Email: shen7@math.purdue.edu. This author's research is partially supported by the U.S. National Science Foundation under grants DMS-1419053, DMS-1620262 and AFOSR grant FA9550-16-1-0102.

where some coefficients of the bulk energy become negative, an extra sixth order polynomial term is needed [6, 7], which couples the three phase variables altogether.

It is challenging to develop efficient schemes while preserving the energy stability to solve the three component Cahn-Hilliard phase-field system, since all three phase variables are nonlinearly coupled. Although a variety of numerical algorithms have been developed for it, most of the existing methods are either first-order accurate in time, or energy unstable, or highly nonlinear, or even the combinations of these features. We refer to [7] for a summary on recent advances about the three-phase models and their numerical approximations. For the developments of numerical methods, we emphasize that the preservation of energy stability laws is critical to capture the correct long time dynamics of the system especially for the preassigned grid size and time step. Furthermore, the energy-law preserving schemes provide flexibility for dealing with stiffness issue in phase-filed models. For the binary phase field model, it is well known that the simple fully implicit or explicit type schemes will induce quite severe stability conditions on the time step so they are not efficient in practice (cf. [17, 18, 43, 44, 52]). Such phenomena appear in the computations of the three component Cahn-Hilliard model as well (cf. the numerical examples in [7]). In particular, the authors of [7] concluded that (i) the fully implicit discretization of the six-order polynomial term leads to convergence of the Newton linearization method only under very tiny time step practically; (ii) it is an open problem about how to prove the existence and convergence of the numerical solutions for the fully implicit type scheme theoretically for the total spreading case; (iii) it is questionable to establish convex-concave decomposition for the six order polynomial term; and (iv) the semi-implicit scheme is the best choice that the authors in [7] can obtain since it is unconditionally energy stable for arbitrary time step, and the existence and the convergence can be thereby proved. However, the semi-implicit schemes referred in [7] are nonlinear, thus they require some efficient iterative solvers in the implementations.

The focus of this paper is to develop stable (unconditionally energy stability) and more efficient (linear) schemes to solve the three component Cahn-Hilliard system. Instead of using traditional discretization approaches like simple implicit, stabilized explicit, convex splitting, or other various tricky Taylor expansions to discretize the nonlinear potentials, we adopt the *Invariant Energy Quadratization* (IEQ) approach, that has been successfully applied in the context of other models in the authors' other work (cf. [23, 54, 58, 59, 62, 70]). However, the application of IEQ method to the three component model is faced with new challenges due to the particular nonlinearities including the Lagrange multiplier term, and the sixth order polynomial potential. The essential idea of the IEQ method is to transform the free energy into a quadratic form of a new variable via a change of variables since the nonlinear potential is usually bounded from below. The new, equivalent system still retains a similar energy dissipation law in terms of new variables. For the time-continuous case, the energy law of the new reformulated system is equivalent to the energy law of the original system. A great advantage of such a reformulation is that all nonlinear terms can be treated semi-explicitly, which in turn leads to a linear system. Moreover, the resulting linear system is symmetric positive definite, thus it can be solved efficiently with simple iterative methods such as CG or other Krylov subspace methods. Using this new strategy, we develop a series of linear and energy stable numerical schemes, without introducing artificial stabilizers as in [43, 52] or using convex splitting approach (cf. [10, 16, 50, 51]).

In summary, the new numerical schemes that we develop in this paper possess the following properties: (i) the schemes are *accurate* (ready for second order or even higher order in time); (ii) they are *unconditionally energy stable*; and (iii) they are *efficient and easy to implement* (lead to symmetric positive definite linear system at each time step). To the best of our knowledge, the proposed schemes are the first such schemes for solving the three-phase Cahn-Hilliard phase-field system that can have all these desired properties. In addition, when the three component Cahn-Hilliard model is coupled with the hydrodynamics (Navier-Stokes), the proposed schemes can be easily applied without any essential difficulties.

The rest of the paper is organized as follows. In Section 2, we give a brief introduction of the three component Cahn-Hilliard model. In Section 3, we construct numerical schemes and prove their unconditional energy stability and symmetric positivity for the linear system in the time discrete case. In Section 4, we present various 2D and 3D numerical simulations to validate the accuracy and efficiency of the proposed schemes. Finally, some concluding remarks are presented in Section 5.

**2. Model System.** We now introduce the three component Cahn-Hilliard phase-field model proposed in [6, 7]. Let $\Omega$ be a smooth, open bounded, connected domain in $\mathbb{R}^d$, $d = 2, 3$. Let $c_i$ ($i = 1, 2, 3$) be the $i-th$ phase function (or order parameter) which represents the volume fraction of the $i-$th component in the fluid mixture, i.e.,

$$c_i = \begin{cases} 1 & \text{inside the i-th component,} \\ 0 & \text{outside the i-th component.} \end{cases} \tag{2.1}$$

In the phase-field framework, a thin (of thickness $\epsilon$) but smooth layer is used to connect the interface that is between 0 and 1. The three unknowns $c_1, c_2, c_3$ are linked though the relationship:

$$c_1 + c_2 + c_3 = 1. \tag{2.2}$$

This is the link condition for the vector $\mathbf{C} = (c_1, c_2, c_3)$, where it belongs to the hyperplane

$$S = \{\mathbf{C} = (c_1, c_2, c_3) \in \mathbb{R}^3, c_1 + c_2 + c_3 = 1\}. \tag{2.3}$$

In the two-phase model, the free energy of the mixture is as follows,

$$E^{diph}(c) = \int_\Omega \Big(\frac{3}{4}\sigma\epsilon|\nabla c|^2 + 12\frac{\sigma}{\epsilon}c^2(1-c)^2\Big)d\boldsymbol{x}, \tag{2.4}$$

where $\sigma$ is the surface tension parameter, the first term contributes to the hydrophilic type (tendency of mixing) of interactions between the materials and the second part, the double well bulk energy term represents the hydrophobic type (tendency of separation) of interactions. As the consequence of the competition between the two types of interactions, the equilibrium configuration will include a diffusive interface with thickness proportional to the parameter $\epsilon$; and, as $\epsilon$ approaches zero, we expect to recover the sharp interface separating the two different materials (cf., for instance, [8, 13, 64]).

There exist several generalizations from the two-phase model to the three-phase model (cf. [6, 7, 29]). In this paper, we adopt below the approach in [7], where the free energy is defined as:

$$E^{triph}(c_1, c_2, c_3) = \int_\Omega \Big(\frac{3}{8}\Sigma_1\epsilon|\nabla c_1|^2 + \frac{3}{8}\Sigma_2\epsilon|\nabla c_2|^2 + \frac{3}{8}\Sigma_3\epsilon|\nabla c_3|^2 + \frac{12}{\epsilon}F(c_1, c_2, c_3)\Big)d\boldsymbol{x}, \tag{2.5}$$

where the coefficient of entropic terms $\Sigma_i$ can be negative for some specific situations.

To be algebraically consistent with the two-phase systems, the three surface tension parameters $\sigma_{12}, \sigma_{13}, \sigma_{23}$ should verify the following conditions:

$$\Sigma_i = \sigma_{ij} + \sigma_{ik} - \sigma_{jk}, i = 1, 2, 3. \tag{2.6}$$

The nonlinear potential $F(c_1, c_2, c_3)$ is:

$$F(c_1, c_2, c_3) = \sigma_{12}c_1^2c_2^2 + \sigma_{13}c_1^2c_3^2 + \sigma_{23}c_2^2c_3^2 + c_1c_2c_3(\Sigma_1c_1 + \Sigma_2c_2 + \Sigma_3c_3) + 3\Lambda c_1^2c_2^2c_3^2. \tag{2.7}$$

Since $c_1, c_2, c_3$ satisfy the hyperplane link condition (2.2), the free energy can be rewritten as

$$F(c_1, c_2, c_3) = F_0(c_1, c_2, c_3) + P(c_1, c_2, c_3), \tag{2.8}$$

where

$$F_0(c_1, c_2, c_3) = \frac{\Sigma_1}{2} c_1^2 (1-c_1)^2 + \frac{\Sigma_2}{2} c_2^2 (1-c_2)^2 + \frac{\Sigma_3}{2} c_3^2 (1-c_3)^2, \tag{2.9}$$
$$P(c_1, c_2, c_3) = 3\Lambda c_1^2 c_2^2 c_3^2,$$

and $\Lambda$ is a non negative constant.

Therefore, the time evolution of $c_i$ is governed by the gradient of the energy $E^{triph}$ with respect to the $H^{-1}(\Omega)$ gradient flow, namely, the Cahn-Hilliard type dynamics as

$$c_{it} = \frac{M_0}{\Sigma_i} \Delta \mu_i, \tag{2.10}$$

$$\mu_i = -\frac{3}{4} \epsilon \Sigma_i \Delta c_i + \frac{12}{\epsilon} \partial_i F + \beta, \ i = 1, 2, 3, \tag{2.11}$$

with the initial condition

$$c_i|_{(t=0)} = c_i^0, \ i = 1, 2, 3, \ c_1^0 + c_2^0 + c_3^0 = 1, \tag{2.12}$$

where $\beta$ is the Lagrange multiplier to ensure the hyperplane link condition (2.2), that can be derived as

$$\beta = -\frac{4\Sigma_T}{\epsilon} \left( \frac{1}{\Sigma_1} \partial_1 F + \frac{1}{\Sigma_2} \partial_2 F + \frac{1}{\Sigma_3} \partial_3 F \right), \tag{2.13}$$

with

$$\frac{3}{\Sigma_T} = \frac{1}{\Sigma_1} + \frac{1}{\Sigma_2} + \frac{1}{\Sigma_3}. \tag{2.14}$$

We consider in this paper either of the two type boundary conditions below:

$$\text{(i) all variables are periodic, or (ii) } \partial_\mathbf{n} c_i|_{\partial\Omega} = \nabla \mu_i \cdot \mathbf{n}|_{\partial\Omega} = 0, \ i = 1, 2, 3, \tag{2.15}$$

where $\mathbf{n}$ is the unit outward normal on the boundary $\partial\Omega$.

It is easily seen that the three chemical potentials $(\mu_1, \mu_2, \mu_3)$ are linked through the relation

$$\frac{\mu_1}{\Sigma_1} + \frac{\mu_2}{\Sigma_2} + \frac{\mu_3}{\Sigma_3} = 0. \tag{2.16}$$

**Remark** 2.1. *In the physical literatures, the coefficient $\Sigma_i$ is called the spreading coefficient of the phase i at the interface between phases j and k. Since $\Sigma_i$ is determined by the surface tensions $\sigma_{i,j}$, it might not be always positive. If $\Sigma_i > 0$, the spreading is said to be "partial", and if $\Sigma_i < 0$, it is called "total".*

The following lemmas hold (cf. [6]):

LEMMA 2.1. *There exists $\underline{\Sigma} > 0$ such that*

$$\Sigma_1 |\xi_1|^2 + \Sigma_1 |\xi_1|^2 + \Sigma_1 |\xi_1|^2 \geq \underline{\Sigma} \left( |\xi_1|^2 + |\xi_2|^2 + |\xi_3|^2 \right), \ \forall \xi_1 + \xi_2 + \xi_3 = 0, \tag{2.17}$$

*if and only if the following condition holds:*

$$\Sigma_1 \Sigma_2 + \Sigma_1 \Sigma_3 + \Sigma_2 \Sigma_3 > 0, \ \Sigma_i + \Sigma_j > 0, \forall i \neq j. \tag{2.18}$$

LEMMA 2.2. *Let $\sigma_{12}, \sigma_{13}$ and $\sigma_{23}$ be three positive numbers and $\Sigma_1, \Sigma_2$ and $\Sigma_3$ defined by (2.6).*

*For any $\Lambda > 0$, the bulk free energy $F(c_1, c_2, c_3)$ defined in (2.8) is bounded from below if $c_1, c_2, c_3$ is on the hyperplane S in 2D. Furthermore, the lower bound only depends on $\Sigma_1, \Sigma_2, \Sigma_3$ and $\Lambda$.*

**Remark** 2.2. *From Lemma 2.1, when (2.18) holds, the summation of the gradient entropy term is bounded from below since $\nabla(c_1 + c_2 + c_3) = 0$, i.e.,*

$$(2.19) \quad \Sigma_1 \|\nabla c_1\|^2 + \Sigma_2 \|\nabla c_2\|^2 + \Sigma_3 \|\nabla c_3\|^2 \geq \underline{\Sigma}(\|\nabla c_1\|^2 + \|\nabla c_2\|^2 + \|\nabla c_3\|^2) \geq 0.$$

**Remark** 2.3. *The bulk part energy $F(c_1, c_2, c_3)$ defined in (2.8) has to be bounded from below in order to form a meaningful physical system. For partial spreading case ($\Sigma_i > 0 \forall i$), one can drop the six order polynomial term by assuming $\Lambda = 0$ since $F_0(c_1, c_2, c_3) \geq 0$ is naturally satisfied. For the total spreading case, $\Lambda$ has to be non zero. Moreover, to ensure the non-negativity for $F$, $\Lambda$ has to be large enough.*

*For 3D case, it is shown in [6] that the bulk energy $F$ is bounded from below when $P(c_1, c_2, c_3)$ takes the following form:*

$$(2.20) \quad P(c_1, c_2, c_3) = 3\Lambda c_1^2 c_2^2 c_3^2 (\phi_\alpha(c_1) + \phi_\alpha(c_2) + \phi_\alpha(c_3))$$

*where $\phi_\alpha(x) = \frac{1}{(1+x^2)^\alpha}$ with $0 < \alpha \leq \frac{8}{17}$.*

*Since (2.9) is commonly used in literatures (cf. [6, 7]), we adopt it as well for convenience. Nonetheless, it will be clear that the numerical schemes we develop in this paper can deal with either (2.9) or (2.20) without any difficulties.*

**Remark** 2.4. *The system (2.10)-(2.11) is equivalent to the following system with two order parameters,*

$$(2.21) \quad \begin{cases} c_{it} = \frac{M_0}{\Sigma_i} \Delta \mu_i, \\ \mu_i = -\frac{3}{4}\epsilon \Sigma_i \Delta c_i + \frac{12}{\epsilon} \partial_i F + \beta, i = 1, 2, \\ c_3 = 1 - c_1 - c_2, \\ \frac{\mu_3}{\Sigma_3} = -(\frac{\mu_1}{\Sigma_1} + \frac{\mu_2}{\Sigma_2}). \end{cases}$$

*We omit the detailed proof since it is quite similar to Theorem 3.1 in section 3.*

The model equation (2.10)-(2.11) follows the dissipative energy law. More precisely, by taking the $L^2$ inner product of (2.10) with $-\mu_i$, and of (2.11) with $c_{it}$, and perform integration by parts, we obtain

$$(2.22) \quad -(c_{it}, \mu_i) = \frac{M_0}{\Sigma_i} \|\nabla \mu_i\|^2,$$

$$(2.23) \quad (\mu_i, c_{it}) = \frac{3}{4}\epsilon \Sigma_i (\nabla c_i, \partial_t \nabla c_i) + \frac{12}{\epsilon}(\partial_i F, c_{it}) + (\beta, c_{it}).$$

Taking summation of the two equalities for $i = 1, 2, 3$, and notice that $(\beta, (c_1+c_2+c_3)_t) = (\beta, (1)_t) = 0$, we obtain the energy dissipative law as

$$(2.24) \quad \frac{d}{dt} E^{triph}(c_1, c_2, c_3) = -M_0 \Big(\frac{1}{\Sigma_1}\|\nabla \mu_1\|^2 + \frac{1}{\Sigma_2}\|\nabla \mu_2\|^2 + \frac{1}{\Sigma_3}\|\nabla \mu_3\|^2\Big).$$

Since $(\mu_1, \mu_2, \mu_3)$ satisfies the condition (2.16), if (2.18) holds, we can derive

$$(2.25) \quad -M_0\Big(\frac{1}{\Sigma_1}\|\nabla \mu_1\|^2 + \frac{1}{\Sigma_2}\|\nabla \mu_2\|^2 + \frac{1}{\Sigma_3}\|\nabla \mu_3\|^2\Big) \leq -M_0 \underline{\Sigma}\Big(\frac{\|\nabla \mu_1\|^2}{\Sigma_1^2} + \frac{\|\nabla \mu_2\|^2}{\Sigma_2^2} + \frac{\|\nabla \mu_3\|^2}{\Sigma_3^2}\Big) \leq 0.$$

**3. Numerical schemes.** We develop in this section several first and second order, unconditionally energy stable and linear numerical schemes for solving the three component phase-field model (2.10)-(2.11). To this end, the main challenges are how to discretize the following three terms including (i) the nonlinear term associated with the double well potential $F_0$, (ii) the six order polynomial term $P$, (iii) the Lagrange multiplier term $\beta$, especially when some $\Sigma_i < 0$ (total spreading).

We notice that for the two-phase Cahn-Hilliard model, the discretization of the nonlinear, cubic polynomial term induced from the double well potential had been well studied (cf. [27, 32, 34, 37, 49]). In summary, there are two commonly used techniques to discretize it in order to preserve the unconditional energy stability. The first is the so-called convex splitting approach [16, 23, 24, 40, 71, 72], where the convex part of the potential is treated implicitly and the concave part is treated explicitly. The convex splitting approach is energy stable, however, it produces nonlinear schemes, thus the implementations are often complicated with potentially high computational costs.

The second technique is the so-called stabilization approach [9, 33, 36, 41–43, 45–48, 52, 53, 55, 57, 63, 66, 69], where the nonlinear term is treated explicitly. In order to preserve the energy law, a linear stabilizing term has to be added, and the magnitude of that term usually depends on the upper bound of the second order derivative of the G-L potential. The stabilized approach leads purely linear schemes, thus it is easy to implement and solve. However, it appears that second order schemes based on the stabilization are only conditionally energy stability [43]. On the other hand, the nonlinear potential may not satisfy the condition required for the stabilization. A feasible remedy is to make some reasonable revisions to the nonlinear potential in order to obtain a finite bound, for example, the quadratic order cut-off functions for the double well potential (cf. [43,52]). Such method is particularly reliable for those models with maximum principle. If the maximum principle does not hold, modified nonlinear potentials may lead to spurious solutions.

For the three component Cahn-Hilliard model system, the above two approaches cannot be easily applied. First of all, even though the convex-concave decomposition can be applied to $F_0$, it is not clear how to deal with the sixth order polynomial term [7]. Second, it is uncertain whether the solution of the system satisfies certain maximum principle so the condition required for stabilization is not satisfied. Third, unconditional energy stable schemes are hardly obtained for both approaches if second order schemes are considered.

We aim to develop numerical schemes that are efficient (well-posed linear system), stable (unconditionally energy stable), and accurate (ready for second order or even higher order in time). To this end, we use the IEQ approach, without worrying about whether the continuous/discrete maximum principle holds or a convexity/concavity splitting exists.

**3.1. Transformed system.** Since $F(c_1, c_2, c_3)$ is always bounded from below from Lemma 2.2 for 2D and Remark 2.3 for 3D. Therefore, we can rewrite the free energy functional $F(c_1, c_2, c_3)$ to the following equivalent form

$$F(c_1, c_2, c_3) = (F(c_1, c_2, c_3) + B) - B, \tag{3.1}$$

where $B$ is some constant to ensure $F(c_1, c_2, c_3) + B > 0$. Furthermore, we define an auxiliary function as follows,

$$U = \sqrt{F(c_1, c_2, c_3) + B}. \tag{3.2}$$

In turn, the total free energy (2.5) can be rewritten as

$$E^{triph}(c_1, c_2, c_3, U) = \int_\Omega \Big(\frac{3}{8}\Sigma_1\epsilon|\nabla c_1|^2 + \frac{3}{8}\Sigma_2\epsilon|\nabla c_2|^2 + \frac{3}{8}\Sigma_3\epsilon|\nabla c_3|^2 + \frac{12}{\epsilon}U^2\Big)d\boldsymbol{x} - \frac{12}{\epsilon}B|\Omega|, \tag{3.3}$$

Thus we could rewrite the system (2.10)-(2.11) to the following equivalent form:

$$(3.4) \quad c_{it} = \frac{M_0}{\Sigma_i}\Delta\mu_i,$$

$$(3.5) \quad \mu_i = -\frac{3}{4}\epsilon\Sigma_i\Delta c_i + \frac{24}{\epsilon}H_iU + \beta, \ i=1,2,3,$$

$$(3.6) \quad U_t = H_1c_{1t} + H_2c_{2t} + H_3c_{3t},$$

where

$$(3.7) \quad H_1 = \frac{\delta U}{\delta c_1} = \frac{1}{2}\frac{\frac{\Sigma_1}{2}(c_1-c_1^2)(1-2c_1)+6\Lambda c_1 c_2^2 c_3^2}{\sqrt{F(c_1,c_2,c_3)+B}},$$

$$(3.8) \quad H_2 = \frac{\delta U}{\delta c_2} = \frac{1}{2}\frac{\frac{\Sigma_2}{2}(c_2-c_2^2)(1-2c_2)+6\Lambda c_1^2 c_2 c_3^2}{\sqrt{F(c_1,c_2,c_3)+B}},$$

$$(3.9) \quad H_3 = \frac{\delta U}{\delta c_3} = \frac{1}{2}\frac{\frac{\Sigma_3}{2}(c_3-c_3^2)(1-2c_3)+6\Lambda c_1^2 c_2^2 c_3}{\sqrt{F(c_1,c_2,c_3)+B}},$$

$$(3.10) \quad \beta = -\frac{8}{\epsilon}\Sigma_T\Big(\frac{1}{\Sigma_1}H_1 + \frac{1}{\Sigma_2}H_2 + \frac{1}{\Sigma_3}H_3\Big)U.$$

The transformed system (3.4)–(3.6) in the variables $c_1, c_2, c_3, U$ form a closed PDE system with the following initial conditions,

$$(3.11) \quad \begin{cases} c_i(t=0) = c_i^0, i=1,2,3, \\ U(t=0) = U^0 = \sqrt{F(c_1^0,c_2^0,c_3^0)+B}. \end{cases}$$

Since the equations (3.6) for the new variables $U$ is only ordinary differential equation with time, the boundary conditions of the new system (3.4)-(3.6) are still (2.15).

**Remark** 3.1. *The system (3.4)-(3.6) is equivalent to the following two order parameter system*

$$(3.12) \quad \begin{cases} c_{it} = \frac{M_0}{\Sigma_i}\Delta\mu_i, \\ \mu_i = -\frac{3}{4}\epsilon\Sigma_i\Delta c_i + \frac{24}{\epsilon}H_iU + \beta, \ i=1,2, \\ U_t = H_1c_{1t} + H_2c_{2t} + H_3c_{3t}, \end{cases}$$

with

$$(3.13) \quad \begin{aligned} c_3 &= 1 - c_1 - c_2, \\ \frac{\mu_3}{\Sigma_3} &= -\Big(\frac{\mu_1}{\Sigma_1} + \frac{\mu_2}{\Sigma_2}\Big). \end{aligned}$$

*Since the proof is quite similar to Theorem 3.1, we omit the proof here.*

We can easily obtain the energy law for the new system (3.4)-(3.6). Taking the $L^2$ inner product of (3.4) with $-\mu_i$, of (3.5) with $\partial_t c_i$, of (3.6) with $-\frac{24}{\epsilon}U$, taking the summation for $i=1,2,3$, and noticing that $(\beta, \partial_t(c_1+c_2+c_3)) = 0$ from Remark 3.1, we still obtain the energy dissipation law as

$$(3.14) \quad \begin{aligned} \frac{d}{dt}E^{triph}(c_1,c_2,c_3,U) &= -M_0\Big(\frac{1}{\Sigma_1}\|\nabla\mu_1\|^2 + \frac{1}{\Sigma_2}\|\nabla\mu_2\|^2 + \frac{1}{\Sigma_3}\|\nabla\mu_3\|^2\Big) \\ &\leq -M_0\underline{\Sigma}\Big(\frac{\|\nabla\mu_1\|^2}{\Sigma_1^2} + \frac{\|\nabla\mu_2\|^2}{\Sigma_2^2} + \frac{\|\nabla\mu_3\|^2}{\Sigma_3^2}\Big) \leq 0. \end{aligned}$$

**Remark** 3.2. *The new transformed system* (3.4)- (3.6) *is equivalent to the original system* (2.10)-(2.11) *since* (3.2) *can be obtained by integrating* (3.6) *with respect to the time. Therefore, the energy law* (3.14) *for the transformed system is exactly the same as the energy law* (2.24) *for the original system.*

*We emphasize that we will develop energy stable numerical schemes for the new transformed system* (3.4)- (3.6). *The proposed schemes will follow a discrete energy law corresponding to* (3.14) *instead of the energy law for the original system* (2.24).

Let $\delta t > 0$ denote the time step size and set $t^n = n\,\delta t$ for $0 \leq n \leq N$ with $T = N\,\delta t$. We also denote by $(f(\boldsymbol{x}), g(\boldsymbol{x})) = \int_\Omega f(\boldsymbol{x})g(\boldsymbol{x})d\boldsymbol{x}$ the $L^2$ inner product of any two functions $f(\boldsymbol{x})$ and $g(\boldsymbol{x})$, and by $\|f\| = \sqrt{(f,f)}$ the $L^2$ norm of the function $f(\boldsymbol{x})$.

**3.2. First order scheme.** We now present the first order time stepping scheme to solve the system (3.4)-(3.6) where the time derivative is discretized based on the first order backward Euler method.

Assuming that $(c_1, c_2, c_3, U)^n$ are already calculated, we compute $(c_1, c_2, c_3, U)^{n+1}$ from the following temporal discrete system:

$$\frac{c_i^{n+1} - c_i^n}{\delta t} = \frac{M_0}{\Sigma_i}\Delta \mu_i^{n+1}, \tag{3.15}$$

$$\mu_i^{n+1} = -\frac{3}{4}\epsilon\Sigma_i\Delta c_i^{n+1} + \frac{24}{\epsilon}H_i^n U^{n+1} + \beta^{n+1},\ i=1,2,3, \tag{3.16}$$

$$U^{n+1} - U^n = H_1^n(c_1^{n+1} - c_1^n) + H_2^n(c_2^{n+1} - c_2^n) + H_3^n(c_3^{n+1} - c_3^n), \tag{3.17}$$

where

$$H_1^n = \frac{1}{2}\frac{\frac{\Sigma_1}{2}(c_1^n - c_1^{n\,2})(1 - 2c_1^n) + 6\Lambda c_1^n c_2^{n\,2} c_3^{n\,2}}{\sqrt{F(c_1^n, c_2^n, c_3^n) + B}}, \tag{3.18}$$

$$H_2^n = \frac{1}{2}\frac{\frac{\Sigma_2}{2}(c_2^n - c_2^{n\,2})(1 - 2c_2^n) + 6\Lambda c_1^{n\,2} c_2^n c_3^{n\,2}}{\sqrt{F(c_1^n, c_2^n, c_3^n) + B}}, \tag{3.19}$$

$$H_3^n = \frac{1}{2}\frac{\frac{\Sigma_3}{2}(c_3^n - c_3^{n\,2})(1 - 2c_3^n) + 6\Lambda c_1^{n\,2} c_2^{n\,2} c_3^n}{\sqrt{F(c_1^n, c_2^n, c_3^n) + B}}, \tag{3.20}$$

$$\beta^{n+1} = -\frac{8}{\epsilon}\Sigma_T\left(\frac{1}{\Sigma_1}H_1^n + \frac{1}{\Sigma_2}H_2^n + \frac{1}{\Sigma_3}H_3^n\right)U^{n+1}. \tag{3.21}$$

The initial conditions are (3.11), and the boundary conditions are

$$\text{(i) all variables are periodic, or (ii) } \partial_{\mathbf{n}} c_i^{n+1}|_{\partial\Omega} = \nabla \mu_i^{n+1} \cdot \mathbf{n}|_{\partial\Omega} = 0,\ i=1,2,3. \tag{3.22}$$

We immediately derive the following result which ensures the numerical solutions satisfy the hyperplane condition (2.2).

THEOREM 3.1. *The system of* (3.15)-(3.17) *is equivalent to the following scheme with two order parameters,*

$$\frac{c_i^{n+1} - c_i^n}{\delta t} = \frac{M_0}{\Sigma_i}\Delta\mu_i^{n+1}, \tag{3.23}$$

$$\mu_i^{n+1} = -\frac{3}{4}\epsilon\Sigma_1\Delta c_i^{n+1} + \frac{24}{\epsilon}H_i^n U^{n+1} + \beta^{n+1},\ i=1,2, \tag{3.24}$$

$$U^{n+1} - U^n = H_1^n(c_1^{n+1} - c_1^n) + H_2^n(c_2^{n+1} - c_2^n) + H_3^n(c_3^{n+1} - c_3^n), \tag{3.25}$$

with

$$c_3^{n+1} = 1 - c_1^{n+1} - c_2^{n+1}, \tag{3.26}$$

$$\frac{\mu_3^{n+1}}{\Sigma_3} = -(\frac{\mu_1^{n+1}}{\Sigma_1} + \frac{\mu_2^{n+1}}{\Sigma_2}). \tag{3.27}$$

*Proof.* From (3.21), we can easily show that the following indentity holds,

$$\frac{24}{\epsilon}(\frac{H_1^n}{\Sigma_1} + \frac{H_2^n}{\Sigma_2} + \frac{H_3^n}{\Sigma_3})U^{n+1} + \beta^{n+1}(\frac{1}{\Sigma_1} + \frac{1}{\Sigma_2} + \frac{1}{\Sigma_3}) = 0. \tag{3.28}$$

- (i) We first assume that (3.23)-(3.27) are satisfied, and derive (3.15)-(3.16). By taking the summation of (3.23) for $i = 1, 2$, and applying (3.26) and (3.27), we obtain

$$\frac{c_3^{n+1} - c_3^n}{\delta t} = \frac{M_0}{\Sigma_3}\Delta\mu_3^{n+1}. \tag{3.29}$$

From (3.28), (3.24), (3.26) and (3.27), we obtain

$$\begin{aligned}
\mu_3^{n+1} &= -\Sigma_3(\frac{\mu_1^{n+1}}{\Sigma_1} + \frac{\mu_2^{n+1}}{\Sigma_2}) \\
&= -\Sigma_3\Big(-\frac{3}{4}\epsilon\Delta c_1^{n+1} - \frac{3}{4}\epsilon\Delta c_2^{n+1} + \frac{24}{\epsilon}(\frac{H_1^n}{\Sigma_1} + \frac{H_2^n}{\Sigma_2})U^{n+1} + \beta^{n+1}(\frac{1}{\Sigma_1} + \frac{1}{\Sigma_2})\Big) \\
&= \frac{3}{4}\epsilon\Sigma_3\Delta c_3^{n+1} + \frac{24}{\epsilon}H_3^n U^{n+1} + \beta^{n+1}.
\end{aligned} \tag{3.30}$$

- (ii) We then assume that (3.15)-(3.16) are satisfied, and derive (3.23)-(3.27). By taking the summation of (3.15) for $i = 1, 2, 3$, we derive

$$\frac{S^{n+1} - S^n}{\delta t} = M_0\Delta\Theta^{n+1}, \tag{3.31}$$

where $S^n = c_1^n + c_2^n + c_3^n$ and $\Theta^{n+1} = \frac{1}{\Sigma_1}\mu_1^{n+1} + \frac{1}{\Sigma_2}\mu_2^{n+1} + \frac{1}{\Sigma_3}\mu_3^{n+1}$. From (3.16) and (3.28), we derive

$$\Theta^{n+1} = -\frac{3}{4}\epsilon\Delta S^{n+1}. \tag{3.32}$$

By taking the $L^2$ inner product of (3.31) with $-\Theta^{n+1}$, of (3.32) with $S^{n+1} - S^n$, and taking the summation of the two equalities above, we derive

$$\frac{3}{8}\epsilon(\|\nabla S^{n+1}\|^2 - \|\nabla S^n\|^2 + \|\nabla S^{n+1} - \nabla S^n\|^2) + \delta t M_0\|\nabla\Theta^{n+1}\|^2 = 0. \tag{3.33}$$

Since $S^n = 1$, the lefthand side of (3.33) is a sum of non negative terms, thus $\nabla S^{n+1} = 0$, and $\nabla\Theta^{n+1} = 0$, i.e., the functions $S^{n+1}$ and $\Theta^{n+1}$ are constants. Then (3.32) leads to $\Theta^{n+1} = 0$, and (3.31) leads to $S^{n+1} = S^n = 1$. Thus we obtain (3.26). By dividing $\Sigma_i$ for (3.16) and taking the summation of it for $i = 1, 2, 3$, and applying (3.28) and (3.26), we obtain (3.27). □

The most impressing part of the above scheme if we notice that the nonlinear coefficient $H_i$ of the new variable $U$ are treated explicitly, which can tremendously speed up the computations in practice. Moreover, we can rewrite the equations (3.17) as follows:

$$U^{n+1} = H_1^n c_1^{n+1} + H_2^n c_2^{n+1} + H_3^n c_3^{n+1} + Q_1^n, \tag{3.34}$$

where $Q_1^n = U^n - H_1^n c_1^n - H_2^n c_2^n - H_3^n c_3^n$. Thus, the system (3.15)-(3.17) can be rewritten as

$$\frac{c_i^{n+1} - c_i^n}{\delta t} = \frac{M_0}{\Sigma_i} \Delta \mu_i^{n+1}, \tag{3.35}$$

$$\mu_i^{n+1} = -\frac{3}{4}\epsilon\Sigma_i \Delta c_i^{n+1} + \frac{24}{\epsilon} H_i^n (H_1^n c_1^{n+1} + H_2^n c_2^{n+1} + H_3^n c_3^{n+1}) \tag{3.36}$$
$$+ \beta^{n+1} + \frac{24}{\epsilon} H_i^n Q_1^n, \; i = 1, 2, 3,$$

THEOREM 3.2. *Assuming* (2.18), *the linear system* (3.35)-(3.36) *for the variable* $\Phi = (c_1^{n+1}, c_2^{n+1}, c_3^{n+1})^T$ *is symmetric (self-adjoint) and positive definite.*

*Proof.* Taking the $L^2$ inner product of (3.35) with 1, we derive

$$\int_\Omega c_i^{n+1} d\boldsymbol{x} = \int_\Omega c_i^n d\boldsymbol{x} = \cdots = \int_\Omega c_i^0 d\boldsymbol{x}. \tag{3.37}$$

Let $\alpha_i^0 = \frac{1}{|\Omega|} \int_\Omega c_i^0 d\boldsymbol{x}$, $\gamma_i^0 = \frac{1}{|\Omega|} \int_\Omega \mu_i^{n+1} d\boldsymbol{x}$, and define

$$\widehat{c}_i^{n+1} = c_i^{n+1} - \alpha_i^0, \; \widehat{\mu}_i^{n+1} = \mu_i^{n+1} - \gamma_i^0. \tag{3.38}$$

Thus, from (3.35)-(3.36), $(\widehat{c}_i^{n+1}, \widehat{\mu}_i^{n+1})$ are the solutions for the following equations with unknowns $(C_i, \mu_i)$,

$$\frac{C_i}{M_0 \delta t} - \frac{1}{\Sigma_i} \Delta \mu_i = \frac{\widehat{c}_i^n}{M_0 \delta t}, \tag{3.39}$$

$$\mu_i + \gamma_i^0 = -\frac{3}{4}\epsilon\Sigma_i \Delta C_i + \frac{24}{\epsilon} H_i^n P_1^{n+1} + \beta^{n+1} + g_i^n, \tag{3.40}$$

where $P_1^{n+1} = H_1^n C_1 + H_2^n C_2 + H_3^n C_3$ and

$$C_1 + C_2 + C_3 = 0, \; \int_\Omega C_i d\boldsymbol{x} = 0, \; \int_\Omega \mu_i d\boldsymbol{x} = 0. \tag{3.41}$$

We define the inverse Laplace operator $u$ (with $\int_\Omega u d\boldsymbol{x} = 0$) $\to v := \Delta^{-1} u$ by

$$\begin{cases} \Delta v = u, \quad \int_\Omega v d\boldsymbol{x} = 0, \\ \text{with the boundary conditions either (i) } v \text{ is periodic, or (ii) } \partial_{\mathbf{n}} v|_{\partial\Omega} = 0 \end{cases} \tag{3.42}$$

Applying $-\Delta^{-1}$ to (3.39) and using (3.40), we obtain

$$-\frac{\Sigma_i}{M_0 \delta t} \Delta^{-1} C_i - \frac{3}{4}\epsilon\Sigma_i \Delta C_i + \frac{24}{\epsilon} H_i^n P_1^{n+1} + \beta^{n+1} - \gamma_i^0 = -\Sigma_i \Delta^{-1} \frac{\widehat{c}_i^n}{M_0 \delta t} - g_i^n, \; i = 1, 2, 3. \tag{3.43}$$

We express the above linear system (3.43) as $\mathbb{A}\mathbf{C} = \mathbf{b}$, where $\mathbf{C} = (C_1, C_2, C_3)^T$.

For any $\Phi = (\phi_1, \phi_2, \phi_3)^T, \Psi = (\psi_1, \psi_2, \psi_3)^T$ with $\sum_{i=1}^3 \phi_i = \sum_{i=1}^3 \psi_i = 0$ and $\int_\Omega \phi_i d\boldsymbol{x} = \int_\Omega \psi_i d\boldsymbol{x} = 0$, we can easily derive

$$(\mathbb{A}\Phi, \Psi) = (\Phi, \mathbb{A}\Psi), \tag{3.44}$$

thus $\mathbb{A}$ is self-adjoint. Meanwhile, we have

$$
\begin{aligned}
(3.45) \quad (\mathbb{A}\Phi, \Phi) = & \frac{1}{M_0 \delta t}(\Sigma_1(-\Delta^{-1}\phi_1, \phi_1) + \Sigma_2(-\Delta^{-1}\phi_2, \phi_2) + \Sigma_3(-\Delta^{-1}\phi_3, \phi_3)) \\
& + \frac{3}{8}\epsilon(\Sigma_1\|\nabla\phi_1\|^2 + \Sigma_2\|\nabla\phi_2\|^2 + \Sigma_3\|\nabla\phi_3\|^2) + \frac{24}{\epsilon}\|H_1^n\phi_1 + H_2^n\phi_2 + H_3^n\phi_3\|^2
\end{aligned}
$$

Let $d_i = \Delta^{-1}\phi_i$, i.e.,

$$
(3.46) \quad \Delta d_i = \phi_i, \quad \int_\Omega d_i d\boldsymbol{x} = 0,
$$

with periodic boundary conditions or $\partial_\mathbf{n} d_i|_{\partial\Omega} = 0$. Therefore, we have

$$
(3.47) \quad (-\Delta^{-1}\phi_i, \phi_i) = \|\nabla d_i\|^2.
$$

Furthermore, $Z = d_1 + d_2 + d_3$ satisfies

$$
(3.48) \quad \Delta Z = 0, \quad \int_\Omega Z d\boldsymbol{x} = 0,
$$

with periodic boundary conditions or $\partial_\mathbf{n} Z|_{\partial\Omega} = 0$. Thus $Z = d_1 + d_2 + d_3 = 0$. From (2.19), we derive

$$
\begin{aligned}
(3.49) \quad (\mathbb{A}\Phi, \Phi) \geq & \frac{1}{M_0 \delta t}\underline{\Sigma}(\|\nabla d_1\|^2 + \|\nabla d_2\|^2 + \|\nabla d_3\|^2) + \frac{3}{8}\underline{\Sigma}\epsilon(\|\nabla\phi_1\|^2 + \|\nabla\phi_2\|^2 + \|\nabla\phi_3\|^2) \\
& + \frac{24}{\epsilon}\|H_1^n\phi_1 + H_2^n\phi_2 + H_3^n\phi_3\|^2 \geq 0,
\end{aligned}
$$

and $(\mathbb{A}\Phi, \Phi) = 0$ if and only if $\Phi = 0$. Thus we conclude the theorem. $\square$

**Remark 3.3.** *We can show the well-posedness of the linear system $\mathbb{A}\mathbf{C} = \mathbf{b}$ from the Lax-Milgram theorem noticing that the linear operator $\mathbb{A}$ is bounded and coercive in $H^1(\Omega)$. We leave the detailed proof to the interested readers since the proof is rather standard.*

The stability result of the first order scheme (3.15)-(3.17) is given below.

**THEOREM 3.3.** *When (2.18) holds, the first order linear scheme (3.15)-(3.17) is unconditionally energy stable, i.e., satisfies the following discrete energy dissipation law:*

$$
(3.50) \quad \frac{1}{\delta t}(E_{1st}^{n+1} - E_{1st}^n) \leq -M_0\underline{\Sigma}\Big(\frac{\|\nabla\mu_1^{n+1}\|^2}{\Sigma_1^2} + \frac{\|\nabla\mu_2^{n+1}\|^2}{\Sigma_2^2} + \frac{\|\nabla\mu_3^{n+1}\|^2}{\Sigma_3^2}\Big).
$$

*where $E_{1st}^n$ is defined by*

$$
(3.51) \quad E_{1st}^n = \frac{3}{8}\Sigma_1\epsilon\|\nabla c_1^n\|^2 + \frac{3}{8}\Sigma_2\epsilon\|\nabla c_2^n\|^2 + \frac{3}{8}\Sigma_3\epsilon\|\nabla c_3^n\|^2 + \frac{12}{\epsilon}\|U^n\|^2 - \frac{12}{\epsilon}B|\Omega|.
$$

*Proof.* Taking the $L^2$ inner product of (3.15) with $-\delta t\mu_i^{n+1}$, we obtain

$$
(3.52) \quad -(c_i^{n+1} - c_i^n, \mu_i^{n+1}) = \delta t\frac{M_0}{\Sigma_i}\|\nabla\mu_i^{n+1}\|^2.
$$

Taking the $L^2$ inner product of (3.16) with $c_i^{n+1} - c_i^n$ and applying the following identities

$$
(3.53) \quad 2(a - b, a) = |a|^2 - |b|^2 + |a - b|^2,
$$

we derive

$$(\mu_i^{n+1}, c_i^{n+1} - c_i^n) = \frac{3}{8}\epsilon\Sigma_i(\|\nabla c_i^{n+1}\|^2 - \|\nabla c_i^n\|^2 + \|\nabla c_i^{n+1} - \nabla c_i^n\|^2)$$
$$+ \frac{24}{\epsilon}(H_i^n U^{n+1}, c_i^{n+1} - c_i^n) + (\beta^{n+1}, c_i^{n+1} - c_i^n). \quad (3.54)$$

Taking the $L^2$ inner product of (3.17) with $\frac{24}{\epsilon} U^{n+1}$, we obtain

$$\frac{12}{\epsilon}\left(\|U^{n+1}\|^2 - \|U^n\|^2 + \|U^{n+1} - U^n\|^2\right)$$
$$= \frac{24}{\epsilon}\left((H_1^n(c_1^{n+1} - c_1^n), U^{n+1}) + (H_2^n(c_2^{n+1} - c_2^n), U^{n+1}) + (H_3^n(c_3^{n+1} - c_3^n), U^{n+1})\right). \quad (3.55)$$

Combining (3.52), (3.54), taking the summation for $i = 1, 2, 3$ and using (3.55) and (3.26), we have

$$\frac{3}{8}\epsilon\sum_{i=1}^{3}\Sigma_i\left(\|\nabla c_i^{n+1}\|^2 - \|\nabla c_i^n\|^2\right) + \frac{3}{8}\epsilon\sum_{i=1}^{3}\Sigma_i\|\nabla c_i^{n+1} - \nabla c_i^n\|^2$$
$$+ \frac{12}{\epsilon}\left(\|U^{n+1}\|^2 - \|U^n\|^2\right) = -M_0(\frac{1}{\Sigma_1}\|\nabla\mu_1^{n+1}\|^2 + \frac{1}{\Sigma_2}\|\nabla\mu_2^{n+1}\|^2 + \frac{1}{\Sigma_3}\|\nabla\mu_3^{n+1}\|^2)$$
$$\leq -M_0\underline{\Sigma}(\|\nabla\mu_1^{n+1}\|^2 + \|\nabla\mu_2^{n+1}\|^2 + \|\nabla\mu_3^{n+1}\|^2) \leq 0. \quad (3.56)$$

Noticing that $\sum_{i=1}^{3}\nabla(c_i^{n+1} - c_i^n) = 0$, we derive

$$\sum_{i=1}^{3}\left(\Sigma_i\|\nabla c_i^{n+1} - \nabla c_i^n\|^2\right) \geq \underline{\Sigma}\sum_{i=1}^{3}\left(\|\nabla c_i^{n+1} - \nabla c_i^n\|^2\right) \geq 0. \quad (3.57)$$

Therefore, the desired result (3.50) is obtained. □

**Remark** 3.4. *We remark that the discrete energy $E_{1st}$ defined in* (3.51) *is bounded from below. Such property is particularly significant. Otherwise, the energy stability does not make any sense.*

*The essential idea of the IEQ method is to transform the complicated nonlinear potentials into a simple quadratic form in terms set of some new variables via a change of variables. Such a simple way of quadratization provides some great advantages. First, the complicated nonlinear potential is transferred to a quadratic polynomial form that is much easier to handle. Second, the derivative of the quadratic polynomial is linear, which provides the fundamental support for linearization method. Third, the quadratic formulation in terms of new variables can automatically keep this property of positivity (or bounded from below) of the nonlinear potentials.*

*We must notice that the convexity is required in the convex splitting approach [16]; and the boundedness for the second order derivative is required in the stabilization approach [43,52]. Compared with those two methods, the IEQ method provides much more flexibilities to treat the complicated nonlinear terms since the only request of it is that the nonlinear potential has to be bounded from below.*

*Meanwhile, the choices of new variables are not unique. For instance, for the partial spreading case where $\Sigma_i > 0, \forall i$, we can define four functions $\widehat{U}_1, \widehat{U}_2, \widehat{U}_3, V$ as follows,*

$$\widehat{U}_i = c_i(1 - c_i), i = 1, 2, 3,$$
$$V = c_1 c_2 c_3. \quad (3.58)$$

Thus the free energy becomes

$$E(\widehat{U}_1, \widehat{U}_2, \widehat{U}_3, V) = \int_\Omega \Big(\frac{3}{8}\Sigma_1\epsilon|\nabla c_1|^2 + \frac{3}{8}\Sigma_2\epsilon|\nabla c_2|^2 + \frac{3}{8}\Sigma_3\epsilon|\nabla c_3|^2 \tag{3.59}$$
$$+ \frac{12}{\epsilon}\Big(\frac{\Sigma_1}{2}\widehat{U}_1^2 + \frac{\Sigma_2}{2}\widehat{U}_2^2 + \frac{\Sigma_3}{2}\widehat{U}_3^2 + 3\Lambda V^2\Big)\Big)d\boldsymbol{x}.$$

For this new definition of the free energy, one can carry about the similar analysis. The details are left to the interested readers.

However, we notice this particular transformation only works for the partial spreading case. For the total spreading case, since for some $i$, $\Sigma_i < 0$, the energy defined in (3.59) may not be bounded from below. This is the particular reason that we define the new variable $U$ in (3.2) in this paper.

As can be seen in this paper, the IEQ approach is able to provide enough flexibilities to derive the equivalent PDE system, and hence the corresponding numerical schemes with unconditional energy stabilities.

**Remark 3.5.** The proposed scheme follows the new energy dissipation law (3.14) instead of the energy law for the originated system (2.24). For time-continuous case, (3.14) and (2.24) are identical. For time-discrete case, the discrete energy law (3.50) is the first order approximation to the new energy law (3.14). Moreover, the discrete energy functional $E_{1st}^{n+1}$ is also the first order approximation to $E(\phi^{n+1})$ (defined in (2.5)), since $U^{n+1}$ is the first order approximations to $\sqrt{F(c_1^{n+1}, c_2^{n+1}, c_3^{n+1}) + B}$, that can be observed from the following facts, heuristically. We rewrite (3.17) as follows,

$$U^{n+1} - \sqrt{F(c_1^{n+1}, c_2^{n+1}, c_3^{n+1}) + B} = U^n - \sqrt{F(c_1^n, c_2^n, c_3^n) + B} + R_{n+1}, \tag{3.60}$$

where $R_{n+1} = O(\sum_{i=1}^3 (c_i^{n+1} - c_i^n)^2)$.

Since $R_k = O(\delta t^2)$ for $0 \le k \le n+1$ and $U^0 = \sqrt{F(c_1^0, c_2^0, c_3^0) + B}$, by mathematical induction we can easily get

$$U^{n+1} = \sqrt{F(c_1^{n+1}, c_2^{n+1}, c_3^{n+1}) + B} + O(\delta t). \tag{3.61}$$

**3.3. Second order scheme based on Crank-Nicolson.** We now present a second order time stepping scheme to solve the system (3.4)-(3.6).

Assuming that $(c_1, c_2, c_3, U)^n$ and $(c_1, c_2, c_3, U)^{n-1}$ are already calculated, we compute $(c_1, c_2, c_3, U)^{n+1}$ from the following temporal discrete system

$$c_{it} = \frac{M_0}{\Sigma_i}\Delta\mu_i^{n+\frac{1}{2}}, \tag{3.62}$$

$$\mu_i^{n+\frac{1}{2}} = -\frac{3}{4}\epsilon\Sigma_i\Delta\frac{c_i^{n+1} + c_i^n}{2} + \frac{24}{\epsilon}H_i^* U^{n+\frac{1}{2}} + \beta^{n+\frac{1}{2}}, i = 1, 2, 3, \tag{3.63}$$

$$U^{n+1} - U^n = H_1^*(c_1^{n+1} - c_1^n) + H_2^*(c_2^{n+1} - c_2^n) + H_3^*(c_3^{n+1} - c_3^n), \tag{3.64}$$

where

$$(3.65) \quad \begin{cases} U^{n+\frac{1}{2}} = \dfrac{U^{n+1} + U^n}{2}, \\ c_i^* = \dfrac{3}{2}c_i^n - \dfrac{1}{2}c_i^{n-1}, \\ H_1^* = \dfrac{1}{2}\dfrac{\frac{\Sigma_1}{2}(c_1^* - c_1^{*2})(1 - 2c_1^*) + 6\Lambda c_1^* c_2^{*2} c_3^{*2}}{\sqrt{F(c_1^*, c_2^*, c_3^*) + B}}, \\ H_2^* = \dfrac{1}{2}\dfrac{\frac{\Sigma_2}{2}(c_2^* - c_2^{*2})(1 - 2c_2^*) + 6\Lambda c_1^{*2} c_2^* c_3^{*2}}{\sqrt{F(c_1^*, c_2^*, c_3^*) + B}}, \\ H_3^* = \dfrac{1}{2}\dfrac{\frac{\Sigma_3}{2}(c_3^* - c_3^{*2})(1 - 2c_3^*) + 6\Lambda c_1^{*2} c_2^{*2} c_3^*}{\sqrt{F(c_1^*, c_2^*, c_3^*) + B}}, \\ \beta^{n+\frac{1}{2}} = -\dfrac{8}{\epsilon}\Sigma_T(\dfrac{1}{\Sigma_1}H_1^* + \dfrac{1}{\Sigma_2}H_2^* + \dfrac{1}{\Sigma_3}H_3^*)U^{n+\frac{1}{2}}. \end{cases}$$

The initial conditions are (3.11), and the boundary conditions are

(3.66) $\quad$ (i) all variables are periodic, or (ii) $\partial_{\mathbf{n}} c_i^{n+1}|_{\partial\Omega} = \nabla \mu_i^{n+\frac{1}{2}} \cdot \mathbf{n}|_{\partial\Omega} = 0$, $i = 1, 2, 3$.

The following theorem ensures the numerical solution $(c_1^{n+1}, c_2^{n+1}, c_3^{n+1})$ always satisfies the hyperplane link condition (2.2).

THEOREM 3.4. *The system* (3.62)-(3.64) *is equivalent to the following scheme with two order parameters,*

$$(3.67) \quad \frac{c_i^{n+1} - c_i^n}{\delta t} = \frac{M_0}{\Sigma_i}\Delta \mu_i^{n+\frac{1}{2}},$$

$$(3.68) \quad \mu_i^{n+\frac{1}{2}} = -\frac{3}{4}\epsilon\Sigma_i\Delta\frac{c_i^{n+1} + c_i^n}{2} + \frac{24}{\epsilon}H_i^* U^{n+\frac{1}{2}} + \beta^{n+\frac{1}{2}}, i = 1, 2,$$

*with*

$$(3.69) \quad c_3^{n+1} = 1 - c_1^{n+1} - c_2^{n+1},$$

$$(3.70) \quad \frac{\mu_3^{n+\frac{1}{2}}}{\Sigma_3} = -(\frac{\mu_1^{n+\frac{1}{2}}}{\Sigma_1} + \frac{\mu_2^{n+\frac{1}{2}}}{\Sigma_2}).$$

*Proof.* The proof is omitted here since it is similar to that for Theorem 3.1. □

Similar to the first orders scheme, we can rewrite the system (3.64) as follows,

$$(3.71) \quad U^{n+1} = H_1^* c_1^{n+1} + H_2^* c_2^{n+1} + H_3^* c_3^{n+1} + Q_2^n,$$

where $Q_2^n = U^n - H_1^* c_1^n - H_2^* c_2^n - H_3^* c_3^n$. Thus, the system (3.62)- (3.64) can be rewritten as

$$(3.72) \quad \frac{c_i^{n+1} - c_i^n}{\delta t} = \frac{M_0}{\Sigma_i}\Delta \mu_i^{n+\frac{1}{2}},$$

$$(3.73) \quad \begin{aligned} \mu_i^{n+\frac{1}{2}} = &-\frac{3}{4}\epsilon\Sigma_i\Delta\frac{c_i^{n+1} + c_i^n}{2} + \frac{24}{\epsilon}H_i^*(H_1^* c_1^{n+1} + H_2^* c_2^{n+1} + H_3^* c_3^{n+1}) \\ &+ \beta^{n+\frac{1}{2}} + \frac{24}{\epsilon}H_i^* Q_2^n, i = 1, 2, 3. \end{aligned}$$

**THEOREM** 3.5. *The linear system* (3.72)-(3.73) *for the variable* $\Phi = (c_1^{n+1}, c_2^{n+1}, c_3^{n+1})^T$ *is self-adjoint and positive definite.*

*Proof.* The proof is omitted here since it is similar to that for Theorem 3.2. □

The stability result of the second order Crank-Nicolson scheme (3.62)-(3.64) is given below.

**THEOREM** 3.6. *Assuming* (2.18), *the second order Crank-Nicolson scheme* (3.62)-(3.64) *is unconditionally energy stable and satisfies the following discrete energy dissipation law:*

$$
\begin{aligned}
\frac{1}{\delta t}(E_{cn}^{n+1} - E_{cn}^n) &= -M_0 \Big(\frac{1}{\Sigma_1}\|\nabla \mu_1^{n+\frac{1}{2}}\|^2 + \frac{1}{\Sigma_2}\|\nabla \mu_2^{n+\frac{1}{2}}\|^2 + \frac{1}{\Sigma_3}\|\nabla \mu_3^{n+\frac{1}{2}}\|^2\Big) \\
&\leq -M_0 \underline{\Sigma} \Big(\frac{\|\nabla \mu_1^{n+\frac{1}{2}}\|^2}{\Sigma_1^2} + \frac{\|\nabla \mu_2^{n+\frac{1}{2}}\|^2}{\Sigma_2^2} + \frac{\|\nabla \mu_3^{n+\frac{1}{2}}\|^2}{\Sigma_3^2}\Big),
\end{aligned}
\tag{3.74}
$$

*where $E_{cn}^n$ that is defined by*

$$
E_{cn}^n = \frac{3}{8}\Sigma_1 \epsilon \|\nabla c_1^n\|^2 + \frac{3}{8}\Sigma_2 \epsilon \|\nabla c_2^n\|^2 + \frac{3}{8}\Sigma_3 \epsilon \|\nabla c_3^n\|^2 + \frac{12}{\epsilon}\|U^n\|^2 - \frac{12}{\epsilon}B|\Omega|.
\tag{3.75}
$$

*Proof.* Taking the $L^2$ inner product of (3.62) with $-\delta t \mu_i^{n+1}$, we obtain

$$
-(c_i^{n+1} - c_i^n, \mu_i^{n+\frac{1}{2}}) = \delta t \frac{M_0}{\Sigma_i}\|\nabla \mu_i^{n+\frac{1}{2}}\|^2.
\tag{3.76}
$$

Taking the $L^2$ inner product of (3.63) with $c_i^{n+1} - c_i^n$, we obtain

$$
\begin{aligned}
(\mu_i^{n+\frac{1}{2}}, c_i^{n+1} - c_i^n) &= \frac{3}{8}\epsilon \Sigma_i (\|\nabla c_i^{n+1}\|^2 - \|\nabla c_i^n\|^2) \\
&\quad + \frac{24}{\epsilon}(H_i^* U^{n+\frac{1}{2}}, c_i^{n+1} - c_i^n) + (\beta^{n+\frac{1}{2}}, c_i^{n+1} - c_i^n).
\end{aligned}
\tag{3.77}
$$

Taking the $L^2$ inner product of (3.17) with $\frac{24}{\epsilon}U^{n+\frac{1}{2}}$, we obtain

$$
\begin{aligned}
\frac{24}{\epsilon}\Big((H_1^*(c_1^{n+1} - c_1^n), U^{n+\frac{1}{2}}) + (H_2^*(c_2^{n+1} - c_2^n), U^{n+\frac{1}{2}}) + (H_3^*(c_3^{n+1} - c_3^n), U^{n+\frac{1}{2}})\Big) \\
= \frac{12}{\epsilon}(\|U^{n+1}\|^2 - \|U^n\|^2)
\end{aligned}
\tag{3.78}
$$

Combining (3.76), (3.77) for $i = 1, 2, 3$ and (3.78), we derive

$$
\begin{aligned}
\frac{3}{8}\epsilon \sum_{i=1}^{3} \Sigma_i \Big(\|\nabla c_i^{n+1}\|^2 - \|\nabla c_i^n\|^2\Big) + \frac{12}{\epsilon}\Big(\|U^{n+1}\|^2 - \|U^n\|^2\Big) \\
= -M_0(\frac{1}{\Sigma_1}\|\nabla \mu_1^{n+\frac{1}{2}}\|^2 + \frac{1}{\Sigma_1}\|\nabla \mu_1^{n+\frac{1}{2}}\|^2 + \frac{1}{\Sigma_1}\|\nabla \mu_1^{n+\frac{1}{2}}\|^2) \\
\leq -M_0 \Sigma_0 (\frac{\|\nabla \mu_1^{n+\frac{1}{2}}\|^2}{\Sigma_1^2} + \frac{\|\nabla \mu_2^{n+\frac{1}{2}}\|^2}{\Sigma_2^2} + \frac{\|\nabla \mu_3^{n+\frac{1}{2}}\|^2}{\Sigma_3^2}).
\end{aligned}
\tag{3.79}
$$

Thus we obtain the desired result (3.74). □

**Remark** 3.6. We remark that the energy law for the second order Crank-Nicolson scheme is "=" instead of "≤" (the first equality in (3.74)).

**3.4. Second order scheme based on BDF.** Now we deveop another type of second order scheme based on the Adam-Bashforth approach (BDF2).

Assumming that $(c_1, c_2, c_3, U)^n$ and $(c_1, c_2, c_3, U)^{n-1}$ are already calculated, we compute $(c_1, c_2, c_3, U)^{n+1}$ from the following discrete system:

$$\text{(3.80)} \qquad \frac{3c_i^{n+1} - 4c_i^n + c_i^{n-1}}{2\delta t} = \frac{M_0}{\Sigma_i} \Delta \mu_i^{n+1},$$

$$\text{(3.81)} \qquad \mu_i^{n+1} = -\frac{3}{4}\epsilon \Sigma_1 \Delta c_i^{n+1} + \frac{24}{\epsilon} H_i^\dagger U^{n+1} + \beta^{n+1}, i = 1, 2, 3,$$

$$\text{(3.82)} \qquad 3U^{n+1} - 4U^n + U^{n-1} = H_1^\dagger(3c_1^{n+1} - 4c_1^n + c_1^{n-1}) \\ + H_2^\dagger(3c_2^{n+1} - 4c_2^n + c_2^{n-1}) + H_3^\dagger(3c_3^{n+1} - 4c_3^n + c_3^{n-1}),$$

where

$$\text{(3.83)} \qquad \begin{cases} c_i^\dagger = 2c_i^n - c_i^{n-1} \\ H_1^\dagger = \frac{1}{2} \frac{\frac{\Sigma_1}{2}(c_1^\dagger - c_1^{\dagger^2})(1 - 2c_1^\dagger) + 6\Lambda c_1^\dagger c_2^{\dagger^2} c_3^{\dagger^2}}{\sqrt{F(c_1^\dagger, c_2^\dagger, c_3^\dagger) + B}} \\ H_2^\dagger = \frac{1}{2} \frac{\frac{\Sigma_2}{2}(c_2^\dagger - c_2^{\dagger^2})(1 - 2c_2^\dagger) + 6\Lambda c_1^{\dagger^2} c_2^\dagger c_3^{\dagger^2}}{\sqrt{F(c_1^\dagger, c_2^\dagger, c_3^\dagger) + B}} \\ H_3^\dagger = \frac{1}{2} \frac{\frac{\Sigma_3}{2}(c_3^\dagger - c_3^{\dagger^2})(1 - 2c_3^\dagger) + 6\Lambda c_1^{\dagger^2} c_2^{\dagger^2} c_3^\dagger}{\sqrt{F(c_1^\dagger, c_2^\dagger, c_3^\dagger) + B}} \\ \beta^{n+1} = -\frac{8}{\epsilon} \Sigma_T (\frac{1}{\Sigma_1} H_1^\dagger + \frac{1}{\Sigma_2} H_2^\dagger + \frac{1}{\Sigma_3} H_3^\dagger) U^{n+1}. \end{cases}$$

The initial conditions are (3.11), and the boundary conditions are

$$\text{(3.84)} \qquad \text{(i) all variables are periodic, or (ii) } \partial_\mathbf{n} c_i^{n+1}|_{\partial \Omega} = \nabla \mu_i^{n+1} \cdot \mathbf{n}|_{\partial \Omega} = 0, i = 1, 2, 3.$$

Similar to the first order scheme and the second order Crank-Nicolson scheme, the hyperplane link condition still holds for this scheme.

THEOREM 3.7. *The system* (3.80)-(3.82) *is equivalent to the following scheme with two order parameters,*

$$\text{(3.85)} \qquad \frac{3c_i^{n+1} - 4c_i^n + c_i^{n-1}}{2\delta t} = \frac{M_0}{\Sigma_i} \Delta \mu_i^{n+1},$$

$$\text{(3.86)} \qquad \mu_i^{n+1} = -\frac{3}{4}\epsilon \Sigma_i \Delta c_i^{n+1} + \frac{24}{\epsilon} H_i^\dagger U^{n+1} + \beta^{n+1}, i = 1, 2,$$

$$\text{(3.87)} \qquad c_3^{n+1} = 1 - c_1^{n+1} - c_2^{n+1},$$

$$\text{(3.88)} \qquad \frac{\mu_3^{n+1}}{\Sigma_3} = -(\frac{\mu_1^{n+1}}{\Sigma_1} + \frac{\mu_2^{n+1}}{\Sigma_2}).$$

*Proof.* The proof is omitted here since it is similar to that for Theorem 3.1. □

Similar to the two previous schemes, we can rewrite the sysytem (3.82) as

$$\text{(3.89)} \qquad U^{n+1} = H_1^\dagger c_1^{n+1} + H_2^\dagger c_2^{n+1} + H_3^\dagger c_3^{n+1} + Q_3^n,$$

where $Q_3^n = U^\ominus - H_1^\dagger c_1^\ominus - H_2^\dagger c_2^\ominus - H_3^\dagger c_3^\ominus$ and $V^\ominus = \frac{4V^n - V^{n-1}}{3}$ for any variable $V$. Thus, the system

(3.80)-(3.82) can be rewritten as

$$\frac{3c_i^{n+1} - 4c_i^n + c_i^{n-1}}{2\delta t} = \frac{M_0}{\Sigma_i}\Delta\mu_i^{n+1}, \tag{3.90}$$

$$\mu_i^{n+1} = -\frac{3}{4}\epsilon\Sigma_i\Delta c_i^{n+1} + \frac{24}{\epsilon}H_i^\dagger(H_1^\dagger c_1^{n+1} + H_2^\dagger c_2^{n+1} + H_3^\dagger c_3^{n+1}) \tag{3.91}$$
$$+\beta^{n+1} + \frac{24}{\epsilon}H_i^\dagger Q_3^n, \ i = 1, 2, 3,$$

THEOREM 3.8. *The linear system* (3.90)-(3.91) *for the variable* $\Phi = (c_1^{n+1}, c_2^{n+1}, c_3^{n+1})^T$ *is self-adjoint and positive definite.*

*Proof.* The proof is omitted here since it is similar to that for Theorem 3.2. □

THEOREM 3.9. *The second order scheme* (3.80)-(3.82) *is unconditionally energy stable and satisfies the following discrete energy dissipation law:*

$$\frac{1}{\delta t}(E_{bdf}^{n+1} - E_{bdf}^n) \leq -M_0\Big(\frac{1}{\Sigma_1}\|\nabla\mu_1^{n+1}\|^2 + \frac{1}{\Sigma_2}\|\nabla\mu_2^{n+1}\|^2 + \frac{1}{\Sigma_3}\|\nabla\mu_3^{n+1}\|^2\Big) \tag{3.92}$$
$$\leq -M_0\underline{\Sigma}\Big(\frac{\|\nabla\mu_1^{n+1}\|^2}{\Sigma_1^2} + \frac{\|\nabla\mu_2^{n+1}\|^2}{\Sigma_2^2} + \frac{\|\nabla\mu_3^{n+1}\|^2}{\Sigma_3^2}\Big),$$

*where* $E_{bdf}^n$ *is defined by*

$$E_{bdf}^n = \frac{3}{8}\Sigma_1\epsilon\Big(\frac{\|\nabla c_1^n\|^2}{2} + \frac{\|2\nabla c_1^n - \nabla c_1^{n-1}\|^2}{2}\Big) + \frac{3}{8}\Sigma_2\epsilon\Big(\frac{\|\nabla c_2^n\|^2}{2} + \frac{\|2\nabla c_2^n - \nabla c_2^{n-1}\|^2}{2}\Big) \tag{3.93}$$
$$+ \frac{3}{8}\Sigma_3\epsilon\Big(\frac{\|\nabla c_3^n\|^2}{2} + \frac{\|2\nabla c_3^n - \nabla c_3^{n-1}\|^2}{2}\Big) + \frac{12}{\epsilon}\Big(\frac{\|U^n\|^2}{2} + \frac{\|2U^n - U^{n-1}\|^2}{2}\Big) - \frac{12}{\epsilon}B|\Omega|.$$

*Proof.* Taking the $L^2$ inner product of (3.80) with $-2\delta t\mu_i^{n+1}$, we obtain

$$-(3c_i^{n+1} - 4c_i^n + c_i^{n-1}, \mu_i^{n+1}) = 2\delta t\frac{M_0}{\Sigma_i}\|\nabla\mu_i^{n+1}\|^2. \tag{3.94}$$

Taking the $L^2$ inner product of (3.80) with $3c_i^{n+1} - 4c_i^n + c_i^{n-1}$, and applying the following identities the following identity

$$2(3a - 4b + c, a) = |a|^2 - |b|^2 + |2a - b|^2 - |2b - c|^2 + |a - 2b + c|^2, \tag{3.95}$$

we derive

$$(\mu_i^{n+1}, 3c_i^{n+1} - 4c_i^n + c_i^{n-1}) = \frac{3}{8}\epsilon\Sigma_i\Big(\|\nabla c_i^{n+1}\|^2 - \|\nabla c_i^n\|^2 + \|2\nabla c_i^{n+1} - \nabla c_i^n\|^2 - \|2\nabla c_i^n - \nabla c_i^{n-1}\|^2\Big)$$
$$+ \frac{3}{8}\epsilon\Sigma_i\|\nabla c_i^{n+1} - 2\nabla c_i^n + \nabla c_i^{n-1}\|^2 \tag{3.96}$$
$$+ \frac{24}{\epsilon}(H_i^\dagger U^{n+1}, 3c_i^{n+1} - 4c_i^n + c_i^{n-1}) + (\beta^{n+1}, 3c_i^{n+1} - 4c_i^n + c_i^{n-1}).$$

Taking the $L^2$ inner product of (3.82) with $\frac{24}{\epsilon}U^{n+1}$, we obtain

$$\begin{aligned}
&\frac{12}{\epsilon}\Big(\|U^{n+1}\|^2 - \|U^n\|^2 + \|2U^{n+1} - U^n\|^2 - \|2U^n - U^{n-1}\|^2 + \|U^{n+1} - 2U^n + U^{n-1}\|^2\Big) \\
&= \frac{24}{\epsilon}\Big((H_1^\dagger(3c_1^{n+1} - 4c_1^n + c_1^{n-1}), U^{n+1}) + (H_2^\dagger(3c_2^{n+1} - 4c_2^n + c_2^{n-1}), U^{n+1}) \\
&\quad + (H_3^\dagger(3c_3^{n+1} - 4c_3^n + c_3^{n-1}), U^{n+1})\Big).
\end{aligned} \quad (3.97)$$

Combining (3.76), (3.77) for $i = 1, 2, 3$ and (3.78), we derive

$$\begin{aligned}
&\frac{3}{8}\epsilon\sum_{i=1}^{3}\Sigma_i\Big(\|\nabla c_i^{n+1}\|^2 + \|2\nabla c_i^{n+1} - \nabla c_i^n\|^2\Big) - \frac{3}{8}\epsilon\sum_{i=1}^{3}\Sigma_i\Big(\|\nabla c_i^n\|^2 + \|2\nabla c_i^n - \nabla c_i^{n-1}\|^2\Big) \\
&\quad + \frac{3}{8}\epsilon\sum_{i=1}^{3}\Sigma_i\Big(\|\nabla c_i^{n+1} - 2\nabla c_i^n + \nabla c_i^{n-1}\|^2\Big) + \frac{12}{\epsilon}\Big(\|U^{n+1}\|^2 + \|2U^{n+1} - U^n\|^2\Big) \\
&\quad - \frac{12}{\epsilon}\Big(\|U^n\|^2 + \|2U^n - U^{n-1}\|^2\Big) + \frac{12}{\epsilon}\|U^{n+1} - 2U^n + U^{n-1}\|^2 \\
&= -2\delta t M_0\Big(\frac{1}{\Sigma_1}\|\nabla\mu_1^{n+1}\|^2 + \frac{1}{\Sigma_2}\|\nabla\mu_2^{n+1}\|^2 + \frac{1}{\Sigma_3}\|\nabla\mu_3^{n+1}\|^2\Big) \\
&\leq -2\delta t M_0 \underline{\Sigma}\Big(\frac{\|\nabla\mu_1^{n+1}\|^2}{\Sigma_1^2} + \frac{\|\nabla\mu_2^{n+1}\|^2}{\Sigma_2^2} + \frac{\|\nabla\mu_3^{n+1}\|^2}{\Sigma_3^2}\Big).
\end{aligned} \quad (3.98)$$

Since $\sum_{i=1}^{3}(\nabla c_i^{n+1} - 2\nabla c_i^n + \nabla c_i^{n-1}) = 0$, from Lemma 2.1, we have

$$\sum_{i=1}^{3}\Big\{\Sigma_i\|\nabla c_i^{n+1} - 2\nabla c_i^n + \nabla c_i^{n-1}\|^2\Big\} \geq \underline{\Sigma}\sum_{i=1}^{3}\Big\{\|\nabla c_i^{n+1} - 2\nabla c_i^n + \nabla c_i^{n-1}\|^2\Big\} \geq 0. \quad (3.99)$$

Therefore, we obtain (3.92) after we drop the unnecessary positive terms in (3.98). □

**Remark** 3.7. *From formal Taylor expansion, we find*

$$\begin{aligned}
&\Big(\frac{\|\phi^{n+1}\|^2 + \|2\phi^{n+1} - \phi^n\|^2}{2\delta t}\Big) - \Big(\frac{\|\phi^n\|^2 + \|2\phi^n - \phi^{n-1}\|^2}{2\delta t}\Big) \\
&\cong \Big(\frac{\|\phi^{n+2}\|^2 - \|\phi^n\|^2}{2\delta t}\Big) + O(\delta t^2) \cong \frac{d}{dt}\|\phi(t^{n+1})\|^2 + O(\delta t^2),
\end{aligned} \quad (3.100)$$

*and*

$$\frac{\|\phi^{n+1} - 2\phi^n + \phi^{n-1}\|^2}{\delta t} \cong O(\delta t^3) \quad (3.101)$$

*for any variable $\phi$. Therefore, the discrete energy law (3.92) is a second order approximation of $\frac{d}{dt}E^{triph}(\phi)$ in (3.14).*

**4. Numerical Simulations.** We present in this section several 2D and 3D numerical experiments using the schemes constructed in the last section. The computational domain is $\Omega = (0,1)^d, d = 2, 3$. We use the second order central finite difference method in space with $128^d$ grid points in all simulations.

Unless otherwise explicitly specified, the default parameters are

$$\eta = 0.03, B = 2, M_0 = 10^{-6}, \Lambda = 7. \quad (4.1)$$

| Coarse $\delta t$ | Fine $\delta t$ | $L^2-$ Error | Order | $L^1-$ Error | Order | $L^\infty-$ Error | Order |
|---|---|---|---|---|---|---|---|
| 0.01 | 0.005 | 7.120e-3 | – | 3.969e-1 | – | 3.950e-4 | – |
| 0.005 | 0.0025 | 3.816e-3 | 0.90 | 2.112e-1 | 0.91 | 2.067e-4 | 0.93 |
| 0.0025 | 0.00125 | 1.999e-3 | 0.93 | 1.101e-1 | 0.94 | 1.062e-4 | 0.96 |
| 0.00125 | 0.000625 | 1.028e-3 | 0.96 | 5.645e-2 | 0.96 | 5.395e-5 | 0.98 |
| 0.000625 | 0.0003125 | 5.225e-4 | 0.98 | 2.863e-2 | 0.98 | 2.719e-5 | 0.99 |
| 0.0003125 | 0.00015625 | 2.635e-4 | 0.99 | 1.442e-2 | 0.99 | 1.365e-5 | 0.99 |
| 0.00015625 | 0.000078125 | 1.323e-4 | 0.99 | 7.239e-3 | 0.99 | 6.842e-6 | 1.0 |

TABLE 4.1
The $L^2, L^1, L^\infty$ numerical errors at $t = 1$ that are computed by the first order scheme LS1 using various temporal resolutions. The order parameters are of (4.1) and $128^2$ grid points are used to discretize the space.

| Coarse $\delta t$ | Fine $\delta t$ | $L^2-$ Error | Order | $L^1-$ Error | Order | $L^\infty-$ Error | Order |
|---|---|---|---|---|---|---|---|
| 0.01 | 0.005 | 3.636e-3 | – | 1.541e-1 | – | 2.460e-4 | – |
| 0.005 | 0.0025 | 1.205e-3 | 1.59 | 4.972e-2 | 1.63 | 9.311e-5 | 1.40 |
| 0.0025 | 0.00125 | 3.590e-4 | 1.75 | 1.453e-2 | 1.78 | 2.999e-5 | 1.63 |
| 0.00125 | 0.000625 | 9.936e-5 | 1.85 | 3.941e-3 | 1.88 | 8.653e-6 | 1.79 |
| 0.000625 | 0.0003125 | 2.626e-5 | 1.92 | 1.030e-3 | 1.93 | 2.336e-6 | 1.89 |
| 0.0003125 | 0.00015625 | 6.759e-6 | 1.96 | 2.636e-4 | 1.97 | 6.080e-7 | 1.94 |
| 0.00015625 | 0.000078125 | 1.715e-7 | 1.98 | 6.671e-5 | 1.98 | 1.551e-7 | 1.97 |

TABLE 4.2
The $L^2, L^1, L^\infty$ numerical errors at $t = 1$ that are computed by the second order scheme LS2-CN using various temporal resolutions. The order parameters are of (4.1) and $128^2$ grid points are used to discretize the space.

**4.1. Accuracy test.** For convenience, we denote the first order scheme (3.15)-(3.17) by LS1, the second order scheme (3.62)-(3.64) by LS2-CN, the scheme (3.80)-(3.82) by LS2-BDF.

We set the initial condition as follows,

$$(4.2) \quad \begin{cases} \phi_3 = \frac{1}{2}\Big(1 + \tanh\Big(\frac{R - 0.15}{\epsilon}\Big)\Big), \\ \phi_1 = \frac{1}{2}(1 - \phi_3)\Big(1 + \tanh\Big(\frac{y - 0.5}{\epsilon}\Big)\Big), \\ \phi_2 = 1 - \phi_1 - \phi_3, \end{cases}$$

| Coarse $\delta t$ | Fine $\delta t$ | $L^2-$ Error | Order | $L^1-$ Error | Order | $L^\infty-$ Error | Order |
|---|---|---|---|---|---|---|---|
| 0.01 | 0.005 | 6.368e-3 | – | 5.374e-1 | – | 5.018e-4 | – |
| 0.005 | 0.0025 | 2.584e-3 | 1.30 | 2.143e-1 | 1.33 | 1.811e-4 | 1.47 |
| 0.0025 | 0.00125 | 9.528e-4 | 1.44 | 7.588e-2 | 1.50 | 5.836e-5 | 1.63 |
| 0.00125 | 0.000625 | 3.172e-4 | 1.59 | 2.443e-2 | 1.64 | 2.244e-5 | 1.38 |
| 0.000625 | 0.0003125 | 9.641e-5 | 1.72 | 7.488e-3 | 1.71 | 7.170e-6 | 1.65 |
| 0.0003125 | 0.00015625 | 2.626e-5 | 1.88 | 2.095e-3 | 1.84 | 1.955e-6 | 1.87 |
| 0.00015625 | 0.000078125 | 7.012e-6 | 1.91 | 6.424e-4 | 1.70 | 5.099e-7 | 1.94 |

TABLE 4.3
*The $L^2, L^1, L^\infty$ numerical errors at $t = 1$ that are computed by the second scheme LS2-BDF using various temporal resolutions. The order parameters are of (4.1) and $128^2$ grid points are used to discretize the space.*

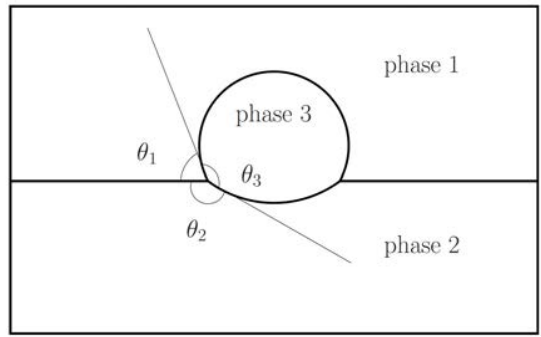

FIG. 4.1. *Theoretical shape of the contact lens at the equilibrium between two stratified fluid components [6].*

where $R = \sqrt{(x - 0.5)^2 + (y - 0.5)^2}$. Since the exact solutions are not known, we compute the errors by adjacent time step. We present the summations of the $L^2, L^1$ and $L^\infty$ errors of the three phase variables at $t = 1$ with different time step sizes in Table 4.1, 4.2 and 4.3 for the three proposed schemes. We observe that the schemes LS1, LS2-CN and LS2-BDF asymptotically match the first order and second order accuracy in time, respectively.

**4.2. Liquid lens between two stratified fluids.** Next, we perform simulations for the classical test case of a liquid lens that is initially spherical sitting at the interface between two other phases. For the accuracy reason, we always use the second order scheme LS2-CN and take the time step $\delta t = 0.001$.

The equilibrium state in the limit $\epsilon \to 0$ can be computed analytically: the final shape of the lens is the union of two pieces of circles, the contact angles being given as a function of the three surface tensions by the Young relations as shown in Fig. 4.1 (cf. [6, 27, 39]).

$$(4.3) \qquad \frac{\sin\theta_1}{\sigma_{23}} = \frac{\sin\theta_2}{\sigma_{13}} = \frac{\sin\theta_3}{\sigma_{12}}.$$

We still use the intitial condition in the previous example, as shown in the first subfigures of Fig. 4.2–4.3, in which, the initial contact angles are $\theta_1 = \theta_2 = \frac{\pi}{2}$ and $\theta_3 = \pi$.

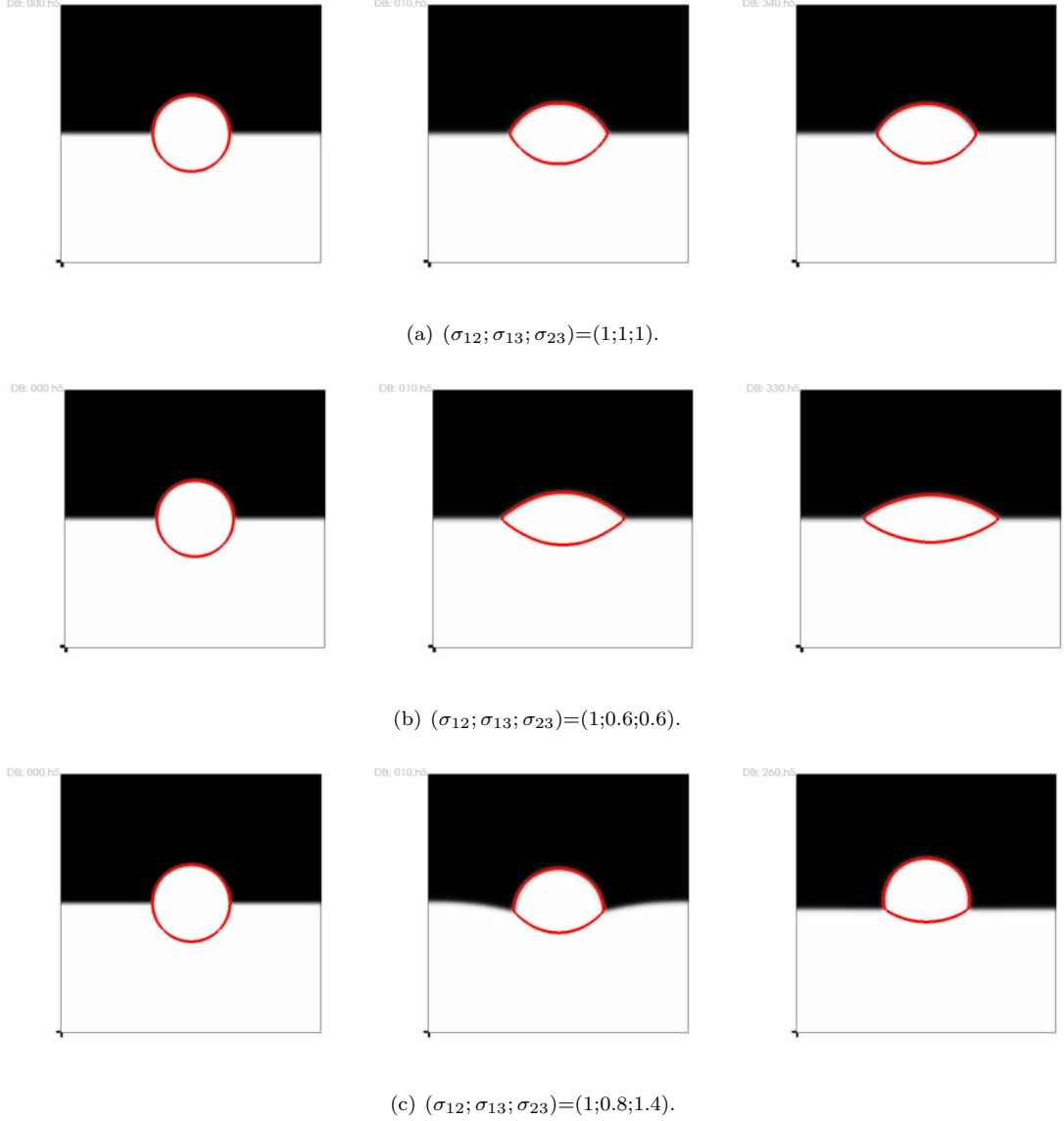

(a) $(\sigma_{12}; \sigma_{13}; \sigma_{23})$=(1;1;1).

(b) $(\sigma_{12}; \sigma_{13}; \sigma_{23})$=(1;0.6;0.6).

(c) $(\sigma_{12}; \sigma_{13}; \sigma_{23})$=(1;0.8;1.4).

FIG. 4.2. *The dynamical behaviors until the steady state of a liquid lens between two stratified fluids for the partial spreading case, with three sets of different surface tension parameters $\sigma_{12}, \sigma_{13}, \sigma_{23}$, where the time step is $\delta t = 0.001$ and $128^2$ grid points are used. The color in black (upper half), white (lower half) and red circle (lens) represent fluids I, II, and III respectively.*

We first simulate the partial spreading case. In Fig. 4.2 (a), we set the three surface tension parameters as $\sigma_{12} = \sigma_{13} = \sigma_{23} = 1$, we observe the three contact angles finally become $\frac{2\pi}{3}$ for all, shown in the final subfigure of Fig. 4.2 because the surface tension force between each phase is the same, which is consistent to the theoretical values of sharp interface from (4.3). In Fig. 4.2 (b), we keep the $\sigma_{12} = 1$ and decrease the other two parameters as $\sigma_{13} = \sigma_{23} = 0.6$. From the contact angle formulation (4.3), we have $\theta_1 = \theta_2 > \theta_3$, which is confirmed by the numerical results shown Fig. 4.2 (b). We further vary the two surface tension parameters $\sigma_{13}$ and $\sigma_{23}$ to be 0.8 and 1.4 respectively while keeping $\sigma_{12} = 1$, in Fig. 4.2 (c). After the intermediate dynamical adjustment, the contact

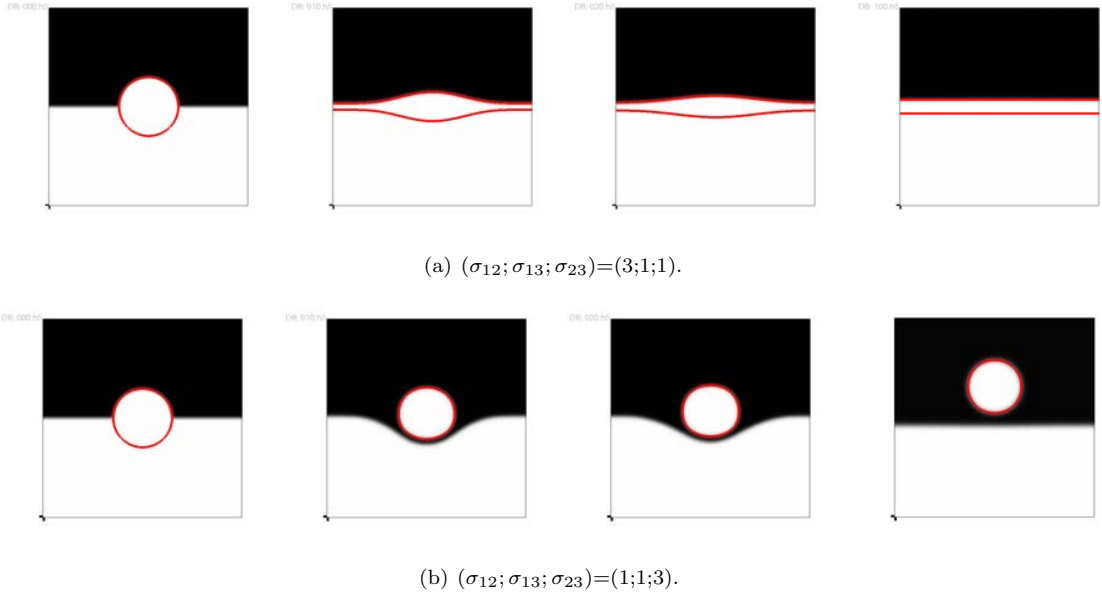

(a) $(\sigma_{12}; \sigma_{13}; \sigma_{23})=(3;1;1)$.

(b) $(\sigma_{12}; \sigma_{13}; \sigma_{23})=(1;1;3)$.

FIG. 4.3. *The dynamical behaviors until the steady state of a liquid lens between two stratified fluids for the total spreading case (no junction points), with two sets of different surface tension parameters $\sigma_{12}, \sigma_{13}, \sigma_{23}$, where the time step is $\delta t = 0.001$ and $128^2$ grid points are used. The color in black (upper half), white (lower half) and red circle (lens) represent fluids I, II, and III respectively.*

angles at equilibrium become $\theta_1 < \theta_3 < \theta_2$, which is consistent to the formulation (4.3) as well.

We then simulate the total spreading case (no junction point) in Fig. 4.3. We set the three surface tension parameters $\sigma_{12} = 3, \sigma_{13} = 1, \sigma_{23} = 1$. From (4.3), the three contact angles can be computed as $\theta_1 = \theta_2 = \pi$ and $\theta_3 = 0$, that can be observed in Fig. 4.3 (a), where the third fluid component $c_3$ totally spread to a layer. For the final case, we set $\sigma_{12} = 1, \sigma_{13} = 1, \sigma_{23} = 3$. By (4.3), it can be computed that the three contact angles at equilibrium are $\theta_1 = 0, \theta_2 = \theta_3 = \pi$, which indicates that the first and third fluid components $c_1, c_3$ are totally spread, and $c_2$ stays inside the first fluid component, as shown in Fig. 4.3 (b). We remark that all numerical results are qualitatively consistent with the computation results obtained in [6, 29].

**4.3. Spinodal decomposition in 2D.** In this example, we study the phase separation behavior, the so-called the spinodal decomposition phenomenon. The process of the phase separation can be studied by considering a homogeneous binary mixture, which is quenched into the unstable part of its miscibility gap. In this case, the spinodal decomposition takes place, which manifests in the spontaneous growth of the concentration fluctuations that leads the system from the homogeneous to the two-phase state. Shortly after the phase separation starts, the domains of the binary components are formed and the interface between the two phases can be specified [4, 12, 73]. For the accuracy reason, we always use the second order scheme LS2-CN and take the time step $\delta t = 0.001$.

The initial condition is taken as the randomly perturbed concentration fields as follows.

$$(4.4) \qquad \phi_i = 0.5 + 0.001 \text{rand}(x,y), \ c_i|_{(t=0)} = \frac{\phi_i}{\phi_1 + \phi_2 + \phi_3}, i = 1, 2, 3,$$

where the $\text{rand}(x,y)$ is the random number in $[-1, 1]$ and has zero mean. To label the three phases, we use pink, grey and yellow to represents phases I, II, and III respectively.

In Fig. 4.4, we perform numerical simulations for the case of order parameters $\sigma_{12} = \sigma_{13} = \sigma_{23} = 1$ as Fig. 4.2 (a). We observe the phase separation behaviors and the final equilibrium solution $t = 30000$ presents a very regular shape where the three contact angles are $\theta_1 = \theta_2 = \theta_3 = \frac{2\pi}{3}$. In Fig. 4.5, with the same initial condition, we set the surface tension parameters are $\sigma_{12} = 1, \sigma_{13} = 0.8, \sigma_{23} = 1.4$. The final equilibrium solution after $t = 30000$ presents show three different contact angles that obey $\theta_1 < \theta_3 < \theta_2$ which is consistent to the example Fig. 4.2 (c). The total spreading case is simulated in Fig. 4.6, in which we set the surface tension parameters are $\sigma_{12} = 1, \sigma_{13} = 1, \sigma_{23} = 3$. The final equilibrium solution after $t = 30000$ presents that no junction points appear, similar as Fig. 4.3 (b).

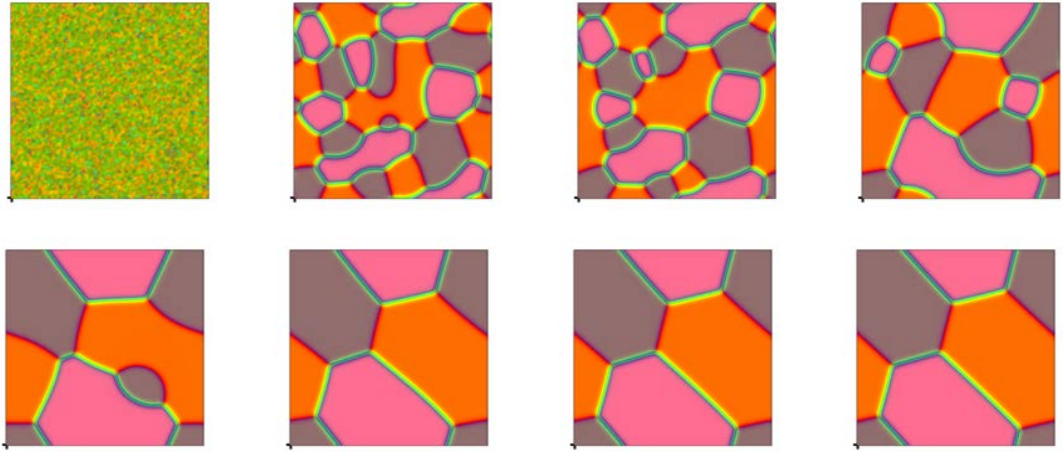

FIG. 4.4. *The 2D dynamical evolution of the three phase variable $c_i, i = 1, 2, 3$ for the partial spreading case, where order parameters are $(\sigma_{12}; \sigma_{13}; \sigma_{23}) = (1 : 1 : 1)$, the time step is $\delta t = 0.001$ and $128^2$ grid points are used. Snapshots of the numerical approximation are taken at $t = 0$, 1000, 2000, 5000, 10000, 20000, 25000, 30000. The color in pink, grey and yellow represents the three phases I, II, and III respectively.*

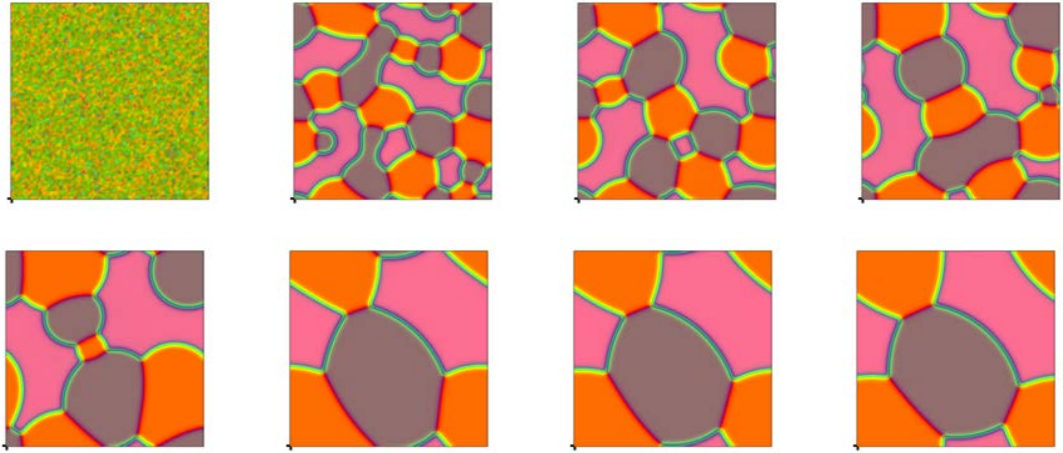

FIG. 4.5. *The 2D dynamical evolution of the three phase variables $c_i, i = 1, 2, 3$ for the partial spreading case, where the order parameters are $(\sigma_{12}; \sigma_{13}; \sigma_{23}) = (1, 0.8, 1.4)$, the time step is $\delta t = 0.001$ and $128^2$ grid points are used. Snapshots of the numerical approximation are taken at $t = 0$, 1000, 2000, 5000, 10000, 20000, 25000, 30000. The color in pink, grey and yellow represents the three phases I, II, and III respectively.*

In Figure 4.7, we present the evolution of the free energy functional for all three cases. The energy curves show the decays with time that confirms that our algorithms are unconditionally stable.

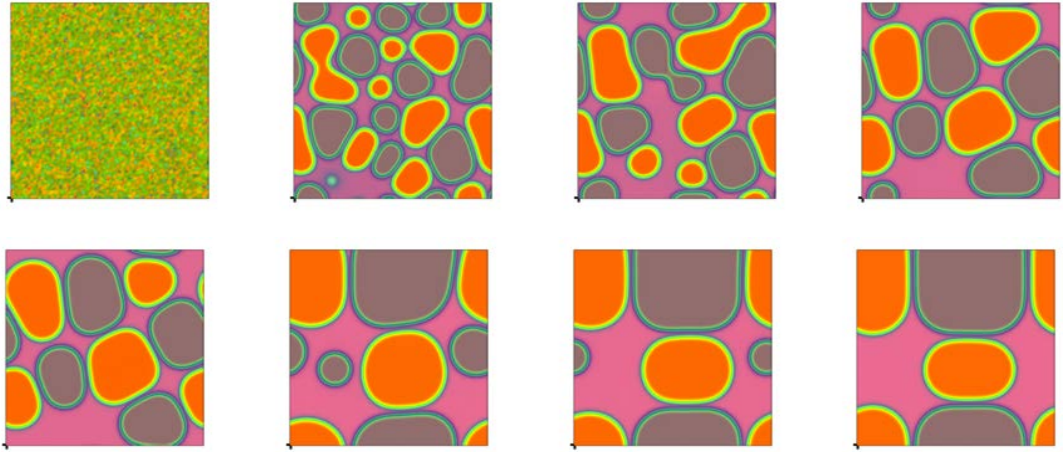

FIG. 4.6. *The 2D dynamical evolution of the three phase variable $c_i$ for the total spreading case (no junction points), where the order parameters are $(\sigma_{12}; \sigma_{13}; \sigma_{23}) = (1, 1, 3)$, the time step is $\delta t = 0.001$ and $128^2$ grid points are used. Snapshots of the numerical approximation are taken at $t = 0$, 1000, 2000, 5000, 10000, 20000, 25000, 30000. The color in pink, grey and yellow represents the three phases I, II, and III respectively.*

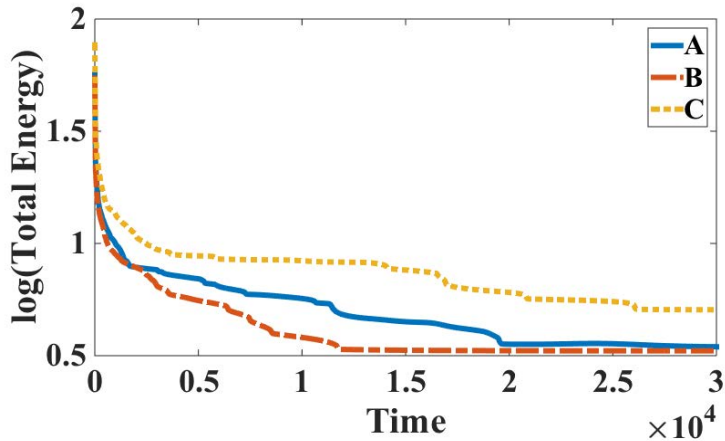

FIG. 4.7. *Time evolution of the free energy functional using the algorithm LS2-CN using $\delta t = 0.001$ for the three order parameters set of $A : (\sigma_{12}; \sigma_{13}; \sigma_{23}) = (1, 0.8, 1.4)$ (partial spreading), $B : (\sigma_{12}; \sigma_{13}; \sigma_{23}) = (1, 1, 1)$ (partial spreading), and $C : (\sigma_{12}; \sigma_{13}; \sigma_{23}) = (1, 1, 3)$ (total spreading). The x-axis is time, and y-axis is $log_{10}$(total energy).*

**4.4. Spinodal decomposition in 3D.** Finally, we present 3D simulations of the phase separation dynamics using the second order scheme LS2-CN and time step $\delta t = 0.001$. In order to be consistent with the 2D case, the initial condition reads as follows,

(4.5)
$$\phi(t = 0) = 0.5 + 0.001 \text{rand}(x, y, z),$$

where the rand$(x, y, z)$ is the random number in $[-1, 1]$ with zero mean.

Fig. 4.8 shows the dynamical behaviors of the phase separation for three equal surface tension parameters $\sigma_{12} = \sigma_{13} = \sigma_{23} = 1$. The dynamical evolutions are consistent with the 2D example Fig.

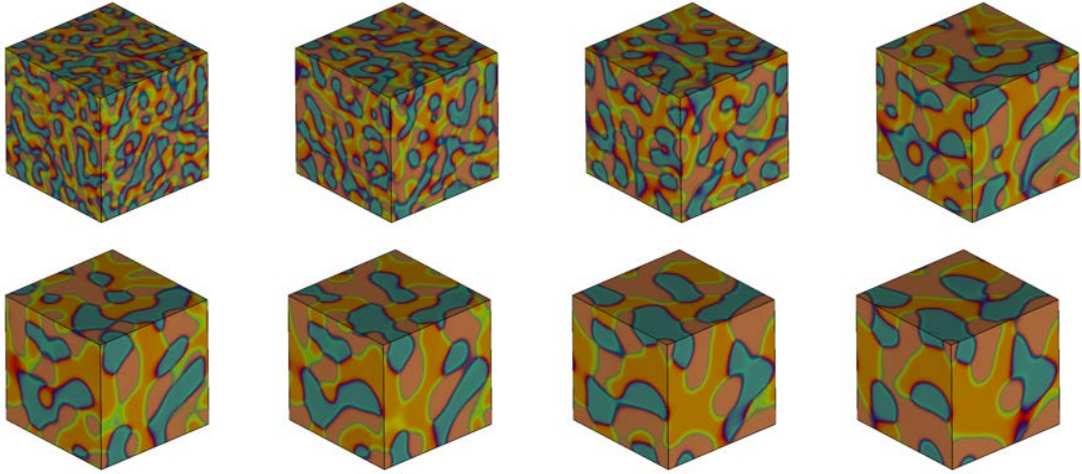

FIG. 4.8. *The 3D dynamical evolution of the three phase variables $c_i, i = 1, 2, 3$ for the partial spreading case, where the order parameters are $(\sigma_{12}; \sigma_{13}; \sigma_{23}) = (1 : 1 : 1)$ and time step is $\delta t = 0.001$. $128^3$ grid points are used to discretize the space. Snapshots of the numerical approximation are taken at $t = 50, 100, 200, 500, 750, 1000, 1500, 2000$. The color in pink, grey and yellow represents the three phases I, II, and III respectively.*

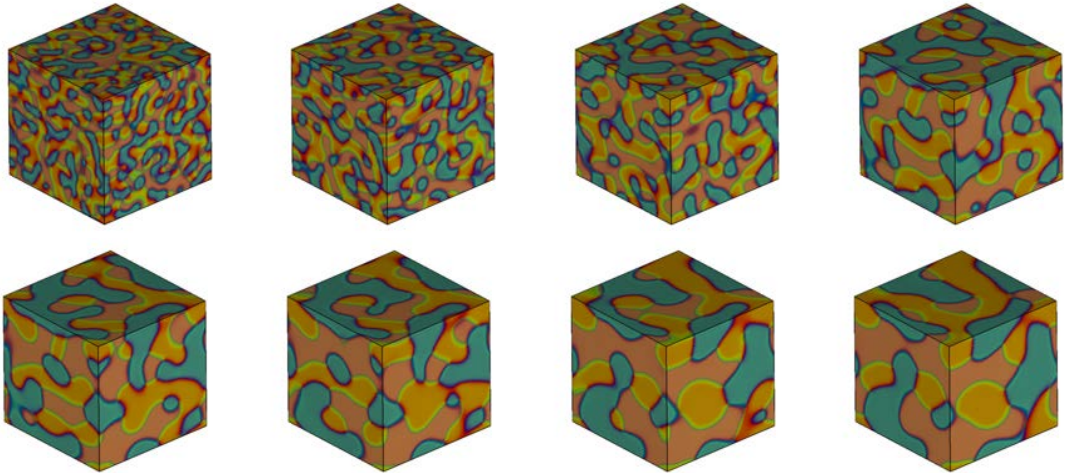

FIG. 4.9. *The 3D dynamical evolution of the three phase variable $c_i, i = 1, 2, 3$ for the partial spreading case, where the order parameters are $(\sigma_{12}; \sigma_{13}; \sigma_{23}) = (1 : 0.8 : 1.4)$ and time step is $\delta t = 0.001$. $128^3$ grid points are used to discretize the space. Snapshots of the numerical approximation are taken at $t = 50, 100, 200, 500, 750, 1000, 1500, 2000$. The color in pink, grey and yellow represents the three phases I, II, and III respectively.*

4.4. We observe that the three components accumulate, as shown in 2D case. In Fig. 4.9, we set the surface tension parameters $\sigma_{12} = 1, \sigma_{13} = 0.8, \sigma_{23} = 1.4$, we observe that the three components accumulate but with different contact angles which is also consistent to the 2D case. In Fig. 4.10, we present the evolution of the free energy functional, in which the energy curves show the decays with time.

**5. Concluding Remarks.** We develop in this paper several efficient time stepping schemes for a three-component Cahn-Hilliard phase-field model that are linear and unconditionally energy stable based on a novel IEQ approach. The proposed schemes bypass the difficulties encountered in the convex splitting and the stabilized approach and enjoy the following desirable properties: (i) *accurate* (up to second order in time); (ii) *unconditionally energy stable*; and (iii) *easy to implement* (one only solves linear equations at each time step). Moreover, the resulting linear system at each time step is symmetric, positive definite so that it can be efficiently solved by any Krylov subspace methods with suitable (e.g., block-diagonal) pre-conditioners.

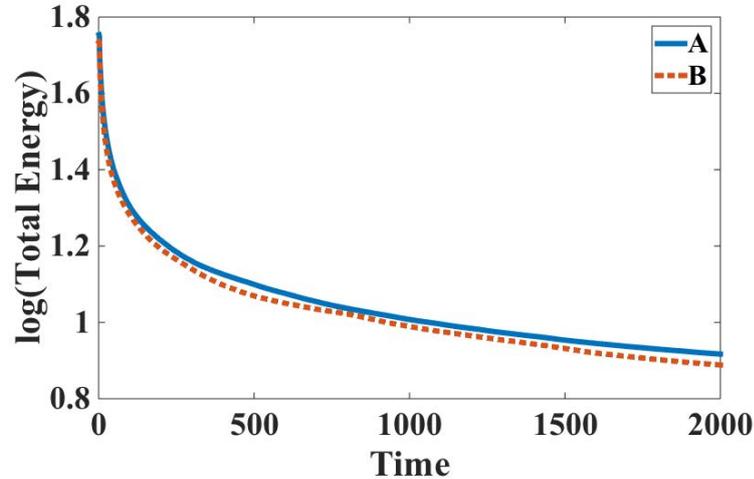

Fig. 4.10. *Time evolution of the free energy functional using the algorithm LS2-CN using $\delta t = 0.001$ for the three order parameters set of $A : (\sigma_{12}; \sigma_{13}; \sigma_{23}) = (1, 0.8, 1.4)$ and $B : (\sigma_{12}; \sigma_{13}; \sigma_{23}) = (1, 1, 1)$. The x-axis is time, and y-axis is $log_{10}$(total energy).*

To the best of our knowledge, these new schemes are the first schemes that are linear and unconditionally energy stable for the three component Cahn-Hilliard phase-field model. These schemes can be applied to the hydro-dynamically coupled three phase model without essential difficulty. Although we considered only time discretization in this study, it is expected that similar results can be established for a large class of consistent finite-dimensional Galerkin approximations since the proofs are all based on a variational formulation with all test functions in the same space as the space of the trial functions.